\documentclass[twoside,11pt]{article}

\usepackage{blindtext}

\usepackage[abbrvbib, preprint]{jmlr2e}

\usepackage[ruled,linesnumbered]{algorithm2e}

\usepackage{lastpage}
\jmlrheading{xx}{2023}{1-\pageref{LastPage}}{8/23}{x/xx}{}{Yusu Hong and Junhong Lin}

\ShortHeadings{High Probability Convergence of Adam}{Hong and Lin}
\firstpageno{1}

\newcommand{\mR}{\mathbb{R}}

\newcommand{\mE}{\mathbb{E}}
\newcommand{\la}{\langle}
\newcommand{\ra}{\rangle}

\newcommand{\ep}{\epsilon}
\newcommand{\gsi}{g_{s,i}}
\newcommand{\bgsi}{\bar{g}_{s,i}}
\newcommand{\xsi}{x_{s,i}}
\newcommand{\bsi}{b_{s,i}}
\newcommand{\asi}{a_{s,i}}
\newcommand{\xisi}{\xi_{s,i}}
\newcommand{\msi}{m_{s,i}}
\newcommand{\hmsi}{\hat{m}_{s,i}}
\newcommand{\vsi}{v_{s,i}}
\begin{document}

\title{High Probability Convergence of Adam Under Unbounded Gradients and Affine Variance Noise}

\author{\name Yusu Hong \email yusuhong@zju.edu.cn \\
       \name Junhong Lin \email junhong@zju.edu.cn \\
       \addr Center for Data Science\\
       Zhejiang University\\
       Hangzhou 310027, P. R. China}
\editor{My editor}

\maketitle
\begin{abstract}
    In this paper, we study the convergence of the Adaptive Moment Estimation (Adam) algorithm under unconstrained non-convex smooth stochastic optimizations. Despite the widespread usage in machine learning areas, its theoretical properties remain limited. Prior researches primarily investigated Adam's convergence from an expectation view, often necessitating strong assumptions like uniformly stochastic bounded gradients or problem-dependent knowledge in prior. As a result, the applicability of these findings in practical real-world scenarios has been constrained. To overcome these limitations, we provide a deep analysis and show that Adam could converge to the stationary point in high probability with a rate of $\mathcal{O}\left({\rm poly}(\log T)/\sqrt{T}\right)$ under coordinate-wise ``affine" variance noise, not requiring any bounded gradient assumption and any problem-dependent knowledge in prior to tune hyper-parameters. Additionally, it is revealed that Adam confines its gradients' magnitudes within an order of $\mathcal{O}\left({\rm poly}(\log T)\right)$. Finally, we also investigate a simplified version of Adam without one of the corrective terms
     and obtain a convergence rate that is adaptive to the noise level. 
\end{abstract}
\section{Introduction} The Adaptive Moment Estimation (Adam) algorithm \citep{kingma2014adam} has witnessed a large success in solving the following unconstrained stochastic optimization problem,
\begin{align*}
    \min_{x \in \mR^d} \mE_{z \in \mathcal{P}}\left[f(x,z)\right],
\end{align*}
where $\mathcal{P}$ is a probability distribution. Similar to most of other adaptive gradient methods, Adam shares the common strength of automatically tuning step-size itself by leveraging historical gradients during the optimization process. Consequently, it alleviates the need for problem-dependent knowledge in prior such as the global Lipschitz constant and the smooth parameter which are often unavailable or difficult to acquire in practical scenarios. Furthermore, owing to its rapid convergence rate and robustness with respect to initial hyper-parameters, Adam and its variant, AdamW, have become the preferred choices in training many tasks in machine learning area, including neural language processing \citep{vaswani2017attention}, computer vision \citep{deng2009imagenet,he2016deep,zhang2020adaptive}, and reinforcement learning \citep{mnih2016asynchronous}.

Despite Adam's widespread usage, its convergence properties exhibit considerable room for improvement especially in non-convex smooth scenarios. For instance, \cite{zou2019sufficient} showed the convergence of Adam when hyper-parameters satisfy proper conditions. \cite{defossez2020simple} and \cite{guo2021novel} established convergence bounds in expectation under assumptions of bounded stochastic gradients and adaptive step-sizes, respectively. However, these assumptions are strict and often not reflective of real-world scenarios. While \cite{zhang2022adam} and \cite{wang2022provable} introduced weaker assumptions, they could only prove that Adam converges in expectation to a bounded region, without guaranteeing convergence to stationary points. Recently, \cite{rakhlinconvergence} established a high probability convergence bound under bounded noise, circumventing the necessity for bounded gradients and extending to the generalized smooth case. Their work revealed that gradients remain bounded along the optimization trajectory only if hyper-parameters are well-tuned. Nonetheless, this convergence result necessitated bounded noise and prior knowledge of smoothness parameter and noise level, limiting its practical applicability in experimental areas. Apart from existing theoretical gaps, there is still a lack of compelling explanations for the selection of two exponential moving average parameters in experiments, where the typical setting is $\beta_1 = 0.9$ and $\beta_2 = 0.999$ \citep{kingma2014adam}.

In this paper, we aim for better understanding the convergence behavior of vanilla Adam that are rarely studied in previous studies and addressing the aforementioned theoretical limitations. We establish a high probability convergence bound for vanilla Adam without requiring any prior knowledge and bounded gradient assumption. We also consider a more general noise model, the coordinate-wise ``affine" variance noise which could cover the commonly used bounded noise and sub-Gaussian noise. Through our analysis, we also show that Adam exhibits sufficient efficacy in ensuring gradients bounded by $\mathcal{O}\left({\rm poly}(\log T)\right)$ throughout the optimization trajectory, leading to convergence to a stationary point with a rate of $\mathcal{O}\left({\rm poly}(\log T)/\sqrt{T}\right)$. This rate is optimal up to a logarithmic factor in the context of non-convex smooth optimization \citep{arjevani2023lower}.

Notably, the parameter setting yielding the optimal rate satisfies $\beta_2 = 1-1/T$ and the learning rate $\eta\sim \mathcal{O}\left(1/\sqrt{T}\right)$, covering all $\beta_1$ being constrained by $0 \le \beta_1 < \beta_2 < 1$ over a total of $T$ iterations, which aligns roughly with the typical setting. Furthermore, the analysis demonstrates that setting $\beta_2 \rightarrow 0$ could potentially lead to divergence, to some extent explaining why a closed to one setting arises in experimental areas. In light of these observations, we bridge the theoretical gap for Adam as highlighted in the second paragraph.

It's worthy noting that we provide a high probability convergence bound with a logarithm dependence on the probability margin. In general, existing literature predominantly focuses on demonstrating convergence from an expectation view, while comparatively few investigations on probabilistic behaviors, perhaps owing to the inherent challenges in theoretical proofs.  However, beyond theoretical interest, the high probability convergence is indispensable for characterizing the convergence behavior of algorithms when they are executed within a limited number of runs, which is a common case in practical settings.

To further pursuit a convergence bound that is adaptive to the noise level, we also study a simplified form of Adam which drops the corrective term for the second moment of stochastic gradient and keeps the other parts unchanged.\footnote{Most previous works studied the simplified Adam without both corrective terms while \cite{defossez2020simple} only dropped the corrective term for the first moment of stochastic gradients. It's worthy noting that our analysis could easily be extended to both of these two simplified forms, potentially leading to better error bounds with simplified analysis.
		We refer to Section \ref{sec:improve} for more discussion.} The theoretical result shows that the simplified form of Adam provides a convergence rate of  $\tilde{\mathcal{O}}\left(1/T + \|\boldsymbol{\sigma}_0\|/\sqrt{T} \right)$ where $\boldsymbol{\sigma}_0 \in \mR^d$ denotes the noise level vector. It thus maintains a property that when the noise level is sufficiently low, the convergence rate accelerates to $\tilde{\mathcal{O}}(1/T)$. Moreover, this convergence rate aligns with the ones observed for vanilla Stochastic Gradient Descent (SGD) in both convex \citep{nemirovski2009robust,agarwal2009information} and non-convex smooth scenarios \citep{ghadimi2013stochastic,arjevani2023lower}.

\subsection{Contributions}
We summarize our main contributions as follows.
\begin{itemize}
    \item We establish a high probability convergence bound for vanilla Adam under non-convex smooth case without requiring bounded gradient assumption and any prior knowledge of problem-dependent parameters. Moreover, we consider the coordinate-wise ``affine" variance noise \citep{wang2023convergence} which covers the bounded noise and sub-Gaussian noise assumption. The result reveals that when exponential moving average parameters satisfy that $\beta_2 = 1-1/T $ and $0\le \beta_1 < \beta_2$, the convergence rate achieves at $\mathcal{O}({\rm poly}(\log T)/\sqrt{T})$ where $T$ denotes the total number of iterations. Both the parameter setting and the attained rate harmonize with previous works \citep{defossez2020simple} where they considered a stronger assumption of bounded stochastic gradients and roughly with the typical setting in experiment.
    \item We also study a simplified version of Adam and provide a high probability convergence bound of the form $\tilde{\mathcal{O}}\left(1/T + \|\boldsymbol{\sigma}_0\|/\sqrt{T} \right)$ that is adaptive to the noise level $\|\boldsymbol{\sigma}_0\|$. 
    \item We show that gradients' magnitudes of Adam along the optimization process are bounded by $\mathcal{O}({\rm poly}(\log T))$ without any prior knowledge to tune the step-size.  
\end{itemize}

The rest of the paper are organized as follows. After introducing some notations, 
in Section \ref{sec:adam}, we introduce the Adam algorithm and some basic assumptions, and then present the high probability convergence results for Adam. 
In Section \ref{sec:relatedwork} we discuss some most relevant literature and comparisons with our results. In Section \ref{sec:proof}, we prove the main convergence result for Adam. In Section \ref{sec:improve}, we establish the convergence bound for the simplified version of Adam. All the missing proofs for some of the technical lemmas are given in the appendix.

\paragraph{Notations}  We use $[d]$ to denote the set $\{1,2,\cdots, d\}$ and $\| \cdot \|, \| \cdot \|_1$ and $\| \cdot \|_{\infty}$ to denote $l_2$-norm, $l_1$-norm and $l_\infty$-norm respectively. $a = \mathcal{O}(b)$ denotes $a \leq C b$ for some universal constant $C$.


\section{High probability convergence for Adam}\label{sec:adam}
We consider unconstrained stochastic optimization over the Euclidean space $\mR^d$ with $l_2$-norm. The objective function $f: \mR^d \rightarrow \mR$ is $L$-smooth satisfying that for any $x, y \in \mR^d$,
\begin{align*}
f(y) - f(x) - \la \nabla f(x), y-x \ra \le \frac{L}{2}\|x-y\|^2.
\end{align*}
Given $x \in \mathbb{R}^d$, we assume a gradient oracle that returns a random vector $g(x)=(g(x)_i)_i \in \mathbb{R}^d$. The deterministic gradient of $f$ at $x$ is denoted by $\nabla f(x) = (\nabla f(x)_i)_i\in \mathbb{R}^d$.

\paragraph{Adam}
We mainly study vanilla Adam that was 1 first put forward by \citep{kingma2014adam}. The algorithm works in the following way:
\begin{algorithm}\label{alg:Adam_scalar}
	\caption{Adam}
	\KwIn{Horizon $T$, $x_1 \in \mathbb{R}^d$, $\beta_1, \beta_2 \in [0,1)$, $m_0 = v_0 = 0$, $\epsilon,\eta > 0$,}
	\For{$s=1,\cdots,T$}{
	Generate $g_s = (g_{s,i})_i = g(x_s) $\;
        \For{$i=1,\cdots, d$}{
        $m_{s,i}=\beta_1 m_{s-1,i} + (1-\beta_1) g_{s,i}$\;
        $v_{s,i} = \beta_2 v_{s-1,i} + (1-\beta_2)g_{s,i}^2 $\;
        $\eta_s = \eta \sqrt{1-\beta_2^s}/(1-\beta_1^s),\ \epsilon_s = \epsilon \sqrt{1-\beta_2^s} $\;
        $x_{s+1,i} = x_{s,i} - \frac{\eta_s}{\sqrt{v_{s,i}}+\epsilon_s} \cdot m_{s,i}$\;
        }
        }
\end{algorithm}

Note that the original Adam algorithm includes two corrective terms for $\msi$ and $v_{s,i}$:
\begin{align}
    \hat{m}_{s,i} = \frac{\msi}{1-\beta_1^s}, \quad \hat{v}_{s,i} = \frac{v_{s,i}}{1-\beta_2^s}. \label{eq:hat_msi}
\end{align}
To simplify the notation, we include the corrective terms into the step size $\eta_s$ satisfying, 
\begin{align}
    \eta_s &=  \frac{\eta \sqrt{1-\beta_2^s}}{1-\beta_1^s} \le \frac{\eta}{1-\beta_1}. \label{eq:eta_s}
\end{align}
Hence, Algorithm \ref{alg:Adam_scalar} performs equivalently to Adam \citep{kingma2014adam}. If we denote the adaptive part of the step-size as 
\begin{align}
     b_0 =0, \quad b_{s,i} := \sqrt{v_{s,i}} + \ep_s, \forall s \ge 1, \forall i \in [d], \label{eq:bsi}
\end{align}
then Algorithm \ref{alg:Adam_scalar} can be written as for any $i \in [d]$, and for any $s \ge 1$,
\begin{align}
    x_{s+1,i} 
    &= x_{s,i} -  \frac{\eta_s(1-\beta_1)}{b_{s,i}} \cdot  g_{s,i} + \beta_1 \cdot \frac{\eta_s b_{s-1,i}}{\eta_{s-1} b_{s,i}}(x_{s,i} - x_{s-1,i}), \label{eq:iteration}
\end{align}
where we let $x_0 = x_1$ and $\eta_0 > 0$ be a positive constant. Compared with the SGD with momentum, \eqref{eq:iteration} reveals that Adam not only performs with the adaptive step-size, but also adaptively tunes its heavy-ball style momentum through $\beta_1$ and $\beta_2$. From this aspect, the adaptive momentum parameter seems mysterious but indeed leads to a convergence through our analysis.

\paragraph{Assumptions} Throughout the paper, we make the following assumptions.
\begin{itemize}
    \item  \textbf{(A1) Bounded below:} The objective function is bounded below, i.e., there exists $f^* > -\infty$ such that $f(x) \ge f^*, \forall x \in \mR^d$;
    \item \textbf{(A2) Unbiased estimator:} The gradient oracle provides an unbiased estimator of $\nabla f(x)$, i.e., $\forall x \in \mR^d$, $\mE\left[g(x) \mid x \right]=\nabla f(x)$;
    \item \textbf{(A3) Coordinate-wise ``affine" variance noise:} We assume a general noise model where for any given $i \in [d]$, $\left(g(x)_{i}-\nabla f(x)_i\right)^2 \le \sigma_{0,i}^2 + \sigma_{1,i}^2 \nabla f(x)_i^2, \text{a.s.}$, $\forall x \in \mR^d$ for some $\boldsymbol{\sigma}_{0} = (\sigma_{0,i})_i \in {\left(\mR^{+}\right)}^d,\boldsymbol{\sigma}_{1} = (\sigma_{1,i})_i \in {\left(\mR^{+}\right)}^d$.
\end{itemize}
The first two assumptions are standard in stochastic optimization. The third assumption introduces a coordinate-wise version of affine variance noise that was previously introduced by \citep{wang2023convergence}. This version permits the noise's magnitude in each coordinate to increase in proportion to the gradient's coordinate. Importantly, this assumption is notably milder compared to the commonly utilized bounded noise or sub-Gaussian noise assumption.\footnote{The sub-Gaussian noise can be proved to be bounded in high-probability.
	Extending the analysis from the bounded noise to the sub-Gaussian noise is not difficult, as noted in \citep{attia2023sgd}.} An instance of a random oracle satisfying coordinate-wise ``affine" variance noise, but not adhering to the bounded noise assumption, is presented in \citep{guo2021novel}. In this case, $g(x)_i = d \cdot \nabla f(x)_i$ for some $i \in [d]$ that is randomly chosen and the rest coordinates of $g(x)$ are zeroes. For further discussion on affine variance noise, interested readers can refer to \citep{bottou2018optimization,guo2021novel,faw2022power,attia2023sgd,wang2023convergence}.

\paragraph{Main convergence result} The following theorem establishes the high probability convergence result for Adam. The detailed proof could be found in Section \ref{sec:proof}.
\begin{theorem}\label{thm:Adam_scalar_1}
    Given $T \ge 1$. Suppose that $\{x_s\}_{s \in [T]}$ is a sequence generated by Algorithm \ref{alg:Adam_scalar}. If Assumptions \textbf{(A1)},\textbf{(A2)} and \textbf{(A3)} hold and 
    the parameters satisfy that 
    \begin{align}
        0 \le \beta_1 < \beta_2 < 1, \quad \eta = C_0\sqrt{1-\beta_2}, \quad \ep = \ep_0\sqrt{1-\beta_2}, \label{eq:parameter_setting}
    \end{align}
    for some constants $C_0,\ep_0 > 0$, then for any given $\delta \in (0,1)$, it holds that with probability at least $1-\delta$,
    \begin{align*}
        \frac{1}{T}\sum_{s=1}^T \|\nabla f(x_s)\|^2 \leq \mathcal{O}\left\{ \frac{G}{TC_0} \cdot \left( \ep_0 + \sqrt{\frac{\|\boldsymbol{\sigma}_0\|^2 + (1+\|\boldsymbol{\sigma}_1\|_{\infty}^2)G}{1-\beta_2}} \right) \right\},
    \end{align*}
    where $\mathcal{F}(T)$ is a polynomial with respect to $T$ 
    and $G$ satisfies the following order,\footnote{The detailed expressions of $\mathcal{F}(T),G$ could be found in \eqref{eq:log_F(t)_1} and \eqref{eq:G_define}.} 
    \begin{align*}
    &\mathcal{F}(T) \sim \mathcal{O}\left\{1 + \frac{\|\boldsymbol{\sigma}_0\|_{\infty}^2 \cdot T + (1+\|\boldsymbol{\sigma}_1\|_{\infty}^2)\left( \|\nabla f(x_1)\|^2 \cdot T + \frac{dC_0^2}{1-\beta_1/\beta_2} \cdot T^3   \right)}{\ep_0^2(1-\beta_2)}  \right\}\\
    &G \sim  \mathcal{O}\left\{f(x_1)-f^* +\frac{C_0(C_0 + \ep_0 + \|\boldsymbol{\sigma}_0\|_{\infty}) }{\beta_2(1-\beta_1)^2(1-\beta_1/\beta_2)}  
    \cdot 
    d \log \left( \frac{\mathcal{F}(T)+dT}{\delta \beta_2^{T}}\right) \right.\\
    &\quad \left.+ \frac{C_0^2(1+\|\boldsymbol{\sigma}_1\|_{\infty}^2)}{\beta_2^2(1-\beta_1)^4(1-\beta_1/\beta_2)^2}\cdot d^2 \log^2 \left( \frac{\mathcal{F}(T)+dT}{\delta \beta_2^{T}}\right)\right\}.
\end{align*}
\end{theorem}
\paragraph{Dependence to $\beta_2$} It's worthy noting that the selection of $\beta_2$ is usually closed to $1$ both in theoretical \citep{zou2019sufficient,defossez2020simple,rakhlinconvergence} and empirical aspect \citep{kingma2014adam,defossez2020simple}. Thus, we could suppose $\beta_2 \ge C$ for some constant $C \in (1/2,1)$ and the dominated order of $G$ with respect to $\beta_2$ is 
\begin{align}
    G \sim \mathcal{O}\left( \log\frac{1}{\beta_2^T} \right). \label{eq:depend_beta2}
\end{align}
An almost the same dependence to $\beta_2$ could be found in \citep{defossez2020simple} under the  bounded stochastic gradients case. We then obtain the following corollary which establishes the optimal rate of $\mathcal{O}\left({\rm poly}(\log T)/\sqrt{T} \right)$ when setting $\beta_2 = 1-1/T$.
\begin{corollary}\label{coro:Adam_scalar_1}
    Given $T \ge 2$. Under the same conditions and settings in Theorem \ref{thm:Adam_scalar_1}, if we set $\beta_2 = 1-1/T$ and let $\beta_1 \in [0,1)$ be a constant satisfying $\beta_1 < \beta_2$, then it holds that with probability at least $1-\delta$,
    \begin{align*}
         \frac{1}{T}\sum_{s=1}^T\|\nabla f(x_s)\|^2  \leq  \mathcal{O}\left\{\frac{1}{T}\cdot \frac{\ep_0 G_1}{C_0} + \frac{G_1}{C_0}\sqrt{\frac{\|\boldsymbol{\sigma}_0\|^2 + (1+\|\boldsymbol{\sigma}_1\|_{\infty}^2)G_1 }{T}} \right\},
    \end{align*}
where $G_1$ is defined by only replacing $\beta_2^T$ with the constant $\mathrm{e}^2$ in $G$ given by Theorem \ref{thm:Adam_scalar_1}.\footnote{Note that since $\beta_2 \ge 1/2$, it's easy to verify that $G_1 \sim \mathcal{O}(\text{poly}(\log T))$ with respect to $T$.}
\end{corollary}
\begin{remark}
The above high probability results indicate that even in the absence of  problem-dependent knowledge in prior and bounded gradient assumptions, Adam achieves the optimal convergence rate with the parameter configuration specified in \eqref{eq:parameter_setting}, along with $\beta_2=1-1/T$. This setup aligns closely with that proposed in \citep{zou2019sufficient,defossez2020simple} and roughly with the commonly used typical values of $\beta_1 = 0.9$ and $\beta_2 = 0.999$. Furthermore, the relationship to $\beta_2$ highlighted in \eqref{eq:depend_beta2} elucidates that an increasing or sufficiently close-to-1 setting is adequate to guarantee Adam's convergence, while an excessively small value might lead to unbounded gradients and subsequent divergence.

\end{remark}

\paragraph{Dependence to $\beta_1$} It's important to observe that the initial outcome in Theorem \ref{thm:Adam_scalar_1} exhibits a high-order reliance on $1/(1-\beta_1)$, causing the bound considerably large as $\beta_1$ approaches 0. However, the subsequent remark validates that this dependence on $1/(1-\beta_1)$ can be reduced by properly adjusting the step-size $\eta$.
\begin{remark}\label{coro:beta1}
Under the same conditions and settings in Theorem \ref{thm:Adam_scalar_1}, if we fix $\beta_1 \in [0,1)$ and set $C_0 = (1-\beta_1)^{3}$, then the convergence rate obtains the order of $\mathcal{O}((1-\beta_1)^{-3})$ with respect to $\beta_1$. It is a little worse than $\mathcal{O}((1-\beta_1)^{-1})$ in \citep{defossez2020simple}. However, \cite{defossez2020simple} dropped the corrective term for $\msi$ and assumed a uniform bound for stochastic gradients, both of which potentially lead to a better dependence.  Using the analysis of this paper for the Adam type algorithm studied in \citep{defossez2020simple}, one can potentially improve the dependence with respect to $1/(1-\beta_1)$,  even with the more general assumptions under this paper.
\end{remark}

\paragraph{Dependence to $\ep$} The convergence bound is proportional to $\ep=\ep_0\sqrt{1-\beta_2}$ ignoring the logarithm term $\mathcal{F}(T)$. In general, our result allows a sufficiently small $\ep$, roughly aligning with the typical one $\ep=10^{-8}$ in experiment.

\section{Discussion and related work}\label{sec:relatedwork}
In the past decade, there has been a substantial surge of scholarly attention directed towards SGD and its adaptive variants. We refer readers to see \citep{ruder2016overview,bottou2018optimization} for an overview. In this section we briefly introduce some most relevant literature to our researches in the sequel. 

\paragraph{Convergence of Adam} Adam was first introduced and proved convergence on online convex optimization by \citep{kingma2014adam}. Later, \cite{reddi2019convergence} pointed out a gap of the proof in \citep{kingma2014adam} and established a counter-example where Adam does not converge. Despite the flaw in the analysis of Adam, it still becomes the default choice optimizer in many machine learning tasks. There arise a number of researches studying the convergence of Adam over non-convex landscapes. \cite{de2018convergence} obtained the convergence requiring gradient bounded and the sign of all stochastic gradients to keep the same, which seems to be hard to check in reality. \cite{zou2019sufficient} provided a convergence bound for Adam under a specific sufficient condition for hyper-parameters and $\beta_2=1-1/T$. \cite{defossez2020simple} established a convergence bound with bounded stochastic gradient assumption and provided some insights to momentum parameters, specifically improving the convergence order to $\mathcal{O}((1-\beta_1)^{-1})$. Recently, \cite{huang2021super} and \cite{guo2021novel} separately established a convergence bound allowing for large $\beta_1$. However, the bounded assumption for adaptive step-sizes, specifically as $C_l \le \frac{1}{\sqrt{v_t}+\ep} \le C_u$, was considered less realistic. In general, these literature all required at least one strong assumption such as uniform bounded (stochastic) gradients. To get rid of this, \cite{shi2020rmsprop} investigated Adam with randomly shuffled setting under affine growth condition (stochastic gradient norm is bounded by an affine function of deterministic gradient norm) instead of bounded gradients, demonstrating convergence to a bounded region when $\beta_2$ is closed to one.  \cite{zhang2022adam} and \cite{wang2022provable} studied Adam under affine growth condition and affine variance noise respectively, but they could only show the convergence to a neighborhood of critical points. 

It's worthy noting that most literature studied variants of Adam ignoring corrective terms. However, as \cite{kingma2014adam} and \cite{defossez2020simple} highlighted through empirical studies, Adam without corrective terms would cause adverse effects on optimization. Besides, the above mentioned works all studied Adam from an expectation view while recently \cite{rakhlinconvergence} presented a high probability convergence bound for vanilla Adam under the generalized smooth and bounded noise case, without any bounded gradient assumptions. However, their bound necessitates rigorous requirements to tune hyper-parameters, specifically $\beta_1 \ge \mathcal{O}\left(1 -  \frac{1}{\sigma^2 G \sqrt{T\log(1/\delta)}} \right)$ and $\eta \sim \mathcal{O}\left( \frac{1-\beta_1}{GL\sqrt{\log(1/\delta)}} \right)$ where $L,\sigma$ denotes generalized smoothness and noise parameters and $G$ is a sufficiently large constant determined by $L,\sigma$. It seems unclear how to check this condition in experiment. In contrast, we provide a high probability convergence result with rather easy-to-check parameter settings that do not rely on any problem-dependent parameters under a weaker noise assumption.

\paragraph{Convergence of other adaptive methods} Several researches have primarily focused on studying the convergence properties of AdaGrad-Norm, the scalar version of AdaGrad \citep{duchi2011adaptive,streeter2010lesson}. Notably, \cite{li2019convergence} established a convergence rate for AdaGrad-Norm with a delayed step-size that remains independent of the current stochastic gradient, albeit requiring knowledge of the smoothness parameter. \cite{ward2020adagrad} derived a similar rate without the need for prior knowledge but still requiring bounded stochastic gradients. To avoid bounded gradient assumption, \cite{faw2022power} obtained the convergence rate with affine variance noise. Later, \cite{wang2023convergence} improved the rate to adapt the noise level with a tighter dependence to $T$. Recently, \cite{faw2023beyond} and \cite{wang2023convergence} relaxed the smooth condition to the generalized smoothness and obtained a similar convergence rate. However, they required generalized smoothness parameters in advance to tune the step-size. The above literature provided convergence bounds in expectation. In parallel, the following works mainly established the convergence in high probability. \cite{kavis2022high} proved that vanilla AdaGrad-Norm converges with a bound adaptive to the sub-Gaussian noise level, requiring uniform bounded gradients. Recently, \cite{alina2023high} got rid of the bounded gradient assumption while achieving the same result under sub-Gaussian noise. Apart from the light-tail noise, \cite{attia2023sgd} studied the affine variance noise and obtained the similar result. 

The element-wise version of AdaGrad was explored by \citep{duchi2011adaptive}. \cite{li2020high} obtained a high probability convergence result for AdaGrad with momentum using a delayed step-size under sub-Gaussian noise. In addition, \cite{wang2023convergence} and \cite{alina2023high} obtained the convergence bounds for AdaGrad under affine variance noise and sub-Gaussian noise respectively. Meanwhile, \cite{shen2023unified} introduced a weighted AdaGrad with unified momentum covering heavy-ball and Nesterov's acceleration. 

Several works focused on other adaptive methods. For example, \cite{zhou2018convergence} studied the convergence of AMSGrad \citep{reddi2019convergence} under sub-Gaussian noise but requiring a more strict bound for gradient summations along the optimization trajectory. \cite{shi2020rmsprop} proved that with proper hyper-parameters, RMSProp \citep{tieleman2012lecture} could converge to a bounded region without bounded gradient assumption. Recently, \cite{zhou2023win} introduced Nesterov-like acceleration into both Adam and AdamW algorithms.

\section{Analysis for vanilla Adam}\label{sec:proof}
In this section we briefly introduce technical innovations of proof with some intuitive insights. In general, we borrow some ideas from previous literature, e.g. \citep{ward2020adagrad,defossez2020simple,kavis2022high,faw2022power,attia2023sgd,alina2023high} and meanwhile, providing a deep analysis framework for Adam. We mainly address the following three key challenges associated with proving convergence in high probability.

 First, there's the unique structure of Adam. Unlike conventional SGD and other adaptive gradient methods like AdaGrad and AMSGrad, which employ a step-size that decreases over time, Adam employs a non-monotonic step-size. This characteristic poses a challenge as many existing analysis frameworks \citep{kavis2022high,wang2023convergence,attia2023sgd} rely on the deceasing step-size property. Moreover, the incorporation of corrective terms for $\msi$ and $\vsi$ introduces additional complexities in demonstrating convergence.

The second major challenge is the potential unbounded gradients throughout the optimization process. Unlike the approach taken in \citep{rakhlinconvergence} under bounded noise case, where they applied a contradiction argument to prove gradient bounded, we use a induction argument to show that Adam can effectively constrain its gradients no greater than $\mathcal{O}({\rm poly} (\log T))$, even without prior knowledge of problem-dependent parameters and under coordinate-wise ``affine" variance noise. This control over gradient magnitudes plays a pivotal role in driving convergence.

Lastly, as emphasized in \eqref{eq:iteration}, Adam operates with an adaptively adjustable momentum parameter. The classical analysis techniques employed for SGD with momentum, as discussed in \citep{ghadimi2015global,mai2020convergence}, cannot be directly applied to Adam due to this distinct characteristic. 

To simplify the notations, we let $\bar{g}_s = (\bar{g}_{s,i})_i = \nabla f(x_s)$ and $\xi_s = (\xi_{s,i})_i = g_s - \bar{g}_s$. Then we start by introducing two auxiliary sequences $\{p_s\}_{s \in [T]}$ and $\{y_s \}_{ s \in [T]}$ as follows, 
\begin{align}
	p_s = \frac{\beta_1}{1-\beta_1}(x_s - x_{s-1}), \quad y_s = p_s + x_s, \quad \forall s\ge 1. \label{eq:define_y_s}
\end{align}

\subsection{Technical lemmas}\label{sec:tech_lemma}
We first state some technical lemmas that are useful in our analysis. Some of the detailed proofs will be given in the appendix.
\begin{lemma}\label{lem:logT_2} 
    Suppose $\{\alpha_s \}_{s \in [T]}$ is a real number sequence. Given $0 < \beta_1 < \beta_2 \le 1$ and $\varepsilon > 0$, we define $\zeta_s = \sum_{j=1}^s \beta_1^{s-j} \alpha_j$, 
    $\gamma_s = \frac{1}{1-\beta_1^s}\sum_{j=1}^s \beta_1^{s-j} \alpha_j$ and  $\theta_s = \sum_{j=1}^s \beta_2^{s-j} \alpha_j^2$, then for any $t \in [T]$,
    \begin{align*}
        &\sum_{s=1}^t \frac{\zeta_s^2}{\varepsilon + \theta_s} \le \frac{1}{(1-\beta_1)(1-\beta_1/\beta_2)}\left(\log \left(1 + \frac{\theta_t}{\varepsilon} \right) - t \log \beta_2 \right), \\
        &\sum_{s=1}^t \frac{\gamma_s^2}{\varepsilon + \theta_s} \le \frac{1}{(1-\beta_1)^2(1-\beta_1/\beta_2)}\left(\log \left(1 + \frac{\theta_t}{\varepsilon} \right) - t \log \beta_2 \right).
    \end{align*}
\end{lemma}
The proof for the first inequality above can be found in \citep{defossez2020simple}, while the proof for the second part will be given in the appendix.
Next, we introduce a concentration inequality for martingale difference sequence that is useful for achieving the high probability bounds, see \citep{li2020high} for the proof.
\begin{lemma} \label{lem:Azuma}
    Suppose $\{Z_s\}_{s \in [T]}$ is a martingale difference sequence with respect to the filtration $\{\mathcal{F}_s\}_{s \in [T]} $. Assume that for each $s \in [T]$, $\sigma_s$ is a $\mathcal{F}_{s-1}$-measurable random variable satisfying
    \begin{align*}
        \mE\left[ \exp(Z_s^2/\sigma_s^2) \mid \mathcal{F}_{s-1} \right] \le \exp(1).
    \end{align*}
    Then for any $\lambda > 0$, and for any $\delta \in (0,1)$, it holds that
    \begin{align*}
        \mathbb{P}\left(\sum_{s=1}^T Z_s > \frac{1}{\lambda}\log {1\over \delta} + \frac{3}{4}\lambda \sum_{s=1}^T \sigma_s^2 \right) \le \delta.
    \end{align*}
\end{lemma}

We next introduce an useful lemma that was classical in the non-convex smooth optimization area and was used in estimating the function value gaps in the analysis of AdaGrad-Norm from \citep{attia2023sgd}. 
\begin{lemma}\label{lem:gradient_delta_s}
	Suppose that $f$ is $L$-smooth, then for any $x\in \mR^d$ we have
	\begin{align*}
	\|\nabla f(x) \|^2 \le 2L(f(x)-f^*).
	\end{align*}
\end{lemma}

\subsection{Algorithm dependent lemmas}\label{sec:alg_lemma}

We now present a rough uniform bound of gradients as follows.
\begin{lemma}\label{lem:OrderT_bound_gradient}
     Let $\{x_s\}_{s \in [T]}$ be given by Algorithm \ref{alg:Adam_scalar} with $0<\beta_1 < \beta_2 \le 1$ and $f$ be $L$-smooth. Then,
    \begin{align*}
        \left|\frac{m_{s-1,i}}{\sqrt{v_{s-1,i}}}\right| \le \sqrt{\frac{(1-\beta_1)(1-\beta_1^{s-1})}{(1-\beta_2)(1-\beta_1/\beta_2)} },\quad \forall s \in[T], \forall i \in[d].
    \end{align*}
    Consequently, if we set $\eta = C_0\sqrt{1-\beta_2}$ for some constant $C_0 > 0$, then 
    \begin{align*}
         \|\nabla f(x_s)\|_{\infty} \le \|\nabla f(x_s)\|  \le \|\nabla f(x_1) \| + LC_{0}\sqrt{\frac{d}{1-\beta_1/\beta_2} }  \cdot s, \quad \forall s\in [T].
    \end{align*}
\end{lemma}
Combining with Lemma \ref{lem:OrderT_bound_gradient} and the smoothness of $f$, we could estimate the norm gap between $\nabla f(x_s) $ and $\nabla f(y_s)$ as follows.
\begin{lemma}\label{lem:gradient_xs_ys}
    Under the same assumptions of Lemma \ref{lem:OrderT_bound_gradient},
    \begin{align*}
        \|\nabla f(x_s) \| \le \|\nabla f(y_s) \| + \frac{LC_0\beta_1\sqrt{d}}{(1-\beta_1)\sqrt{1-\beta_1/\beta_2}}:=\|\nabla f(y_s) \|+M.
    \end{align*}
\end{lemma}
We also present the following lemma which provides an upper bound of order $\mathcal{O}({\rm poly}(\log T))$ for four summations related to the stochastic gradients and the adaptive step-sizes.
\begin{lemma}\label{lem:sum_1}
Let $\gsi,\msi$ be given by Algorithm \ref{alg:Adam_scalar} and $\hmsi,\bsi$ be defined in \eqref{eq:hat_msi} and \eqref{eq:bsi} with $0 \le \beta_1 < \beta_2 < 1$, then for any $t \in [T]$,
    \begin{align*}
        &\sum_{i=1}^d \sum_{s=1}^t\frac{\gsi^2}{\bsi^2}  \le \frac{d}{1-\beta_2}\log\left( \frac{\mathcal{F}(T)}{\beta_2^T} \right),\\
        &\sum_{i=1}^d \sum_{s=1}^{t-1}\frac{\msi^2}{\bsi^2} \le \frac{d(1-\beta_1)}{(1-\beta_2)(1-\beta_1/\beta_2)}\log\left( \frac{\mathcal{F}(T)}{\beta_2^T} \right), \\
        &\sum_{i=1}^d \sum_{s=1}^{t-1}\frac{m_{s,i}^2}{b_{s+1,i}^2} \le \frac{d(1-\beta_1)}{\beta_2(1-\beta_2)(1-\beta_1/\beta_2)}\log\left( \frac{\mathcal{F}(T)}{\beta_2^T} \right), \\
        &\sum_{i=1}^d \sum_{s=1}^{t-1} \frac{\hmsi^2}{\bsi^2} \le \frac{d }{(1-\beta_2)(1-\beta_1/\beta_2)}\log\left( \frac{\mathcal{F}(T)}{\beta_2^T} \right).
    \end{align*}
    where $\mathcal{F}(T) \sim O(T^3)$ with the explicit expression given in \eqref{eq:log_F(t)_1}.
\end{lemma}

\subsection{Proxy step-size}
Before proving the main result, we introduce the proxy step-size $\asi$ that helps to break the correlation between $\bsi$ and $\gsi$. 
We adopt this idea from \citep{ward2020adagrad} that was initially employed for AdaGrad-Norm under bounded noise case, and later under affine variance noise case \citep{faw2022power,attia2023sgd}. Specifically, $\asi$ is defined as 
\begin{align}
    \asi = \sqrt{\tilde{v}_{s,i}} + \ep_s = \sqrt{\beta_2 v_{s-1,i} + (1-\beta_2)(\sigma_{0,i}^2 + (1+\sigma_{1,i}^2)\bgsi^2)} + \ep_s, \label{eq:prox_as}
\end{align}
Unlike the setting for AdaGrad-Norm where a uniform $a_s$ was employed, it's necessary to introduce a proxy step-size that is different on each coordinate. Moreover, the exponential moving average has also been adopted into the proxy step-size.
Due to the dedicated construction, the following lemma shows that $\asi$ well estimates both $\bsi$ and $b_{s-1,i}$ and the gap is controlled by the noise level and the gradient magnitude, which plays a vital role in estimating the effect of adaptive heavy-ball momentum.
\begin{lemma}\label{lem:gap_as_bs}
    Following the definitions of $\bsi,\asi$ in \eqref{eq:bsi} and \eqref{eq:prox_as} respectively together with Assumption (A3), it holds that for any $i \in [d]$ and $s \ge 1$,
    \begin{align*}
        \left| \frac{1}{\asi} - \frac{1}{\bsi} \right| \le \frac{2\sqrt{1-\beta_2}}{\asi\bsi}\sqrt{\sigma_{0,i}^2 + \sigma_{1,i}^2 \bgsi^2},
    \end{align*}
    and 
    \begin{align*}
        \left| \frac{1}{\asi} - \frac{1}{b_{s-1,i}} \right| \le \frac{\sqrt{1-\beta_2}}{b_{s-1,i}\asi}\cdot \max\left\{ \sqrt{\sigma_{0,i}^2 + (1+\sigma_{1,i}^2)\bgsi^2} + \ep,  \sqrt{v_{s-1,i}} + \ep \right\}.
    \end{align*}
\end{lemma}
\subsection{Bounding gradients}
We obtain the following proposition which states that with only smoothness condition and the coordinate-wise ``affine" variance noise, the gradient magnitude is well controlled with a $\mathcal{O}(\text{poly}(\log T))$ order in high probability. 

\begin{proposition}\label{pro:delta_s_1}
    Under the same conditions in Theorem \ref{thm:Adam_scalar_1}, for any given $\delta \in (0,1)$, it holds that with probability at least $1-\delta$,
    \begin{align}
        \|\nabla f(x_{t+1})\|^2 \le G - L \sum_{i=1}^d\sum_{s=1}^t  \eta_s \cdot \frac{\nabla f(x_s)_i}{\asi}, \quad \forall t \in [T], \label{eq:delta_s_final_2}
    \end{align}
    where $G$ is given in Theorem \ref{thm:Adam_scalar_1} with the explicit expression given in \eqref{eq:G_define}. As a consequence, we obtain that $\|\nabla f(x_t)\|^2 \le G, \forall t \in [T+1]$ with probability at least $1-\delta$.
\end{proposition}

\paragraph{Proof of Proposition \ref{pro:delta_s_1}}
	Following from \eqref{eq:define_y_s} and recalling Algorithm \ref{alg:Adam_scalar},
we first reveal the iteration relationship of $y_s$ as
\begin{align*}
    y_{s+1,i} = y_{s,i} - \eta_s \cdot \frac{\gsi}{\bsi}  +  \frac{\beta_1}{1-\beta_1}\cdot \left( \frac{\eta_s b_{s-1,i}}{\eta_{s-1} \bsi} -1 \right)\cdot(\xsi - x_{s-1,i}),\quad \forall s \ge 1.
\end{align*}
Applying the descent lemma of smoothness and then making a decomposition, we have for any $s \ge 1$,
\begin{align}
    f(y_{s+1}) &\le f(y_s) + \langle \nabla f(y_s), y_{s+1}-y_s \rangle + \frac{L}{2}\|y_{s+1}-y_s\|^2 \nonumber \\
               &\le f(y_s) + \underbrace{\sum_{i=1}^d -\eta_s \cdot \nabla f(y_s)_i \cdot \frac{\gsi}{\bsi}}_{\textbf{A}} + \underbrace{\sum_{i=1}^d \frac{\beta_1}{1-\beta_1}\cdot \nabla f(y_s)_i \cdot \left(\frac{\eta_s b_{s-1,i}}{\eta_{s-1} b_{s,i}} -1 \right)\cdot(x_{s,i} - x_{s-1,i})}_{\textbf{B}}  \nonumber \\
               &\quad+ \underbrace{\frac{L}{2}\sum_{i=1}^d \left(\eta_s \cdot \frac{\gsi}{\bsi}  -  \frac{\beta_1}{1-\beta_1}\cdot \left(  \frac{\eta_{s} b_{s-1,i}}{\eta_{s-1} \bsi} -1 \right)\cdot(\xsi - x_{s-1,i}) \right)^2}_{\textbf{C}}.  \label{eq:A+B+C}
\end{align}
To make the notation short, we
denote
\begin{align}
G_t = \max_{s \in [t]} \|\bar{g}_s\|^2,\quad 
\bar{\mathcal{G}}_s = \sqrt{\|\boldsymbol{\sigma}_0\|_{\infty}^2 + \|\boldsymbol{\sigma}_1\|^2\|\bar{g}_s\|_\infty^2 }, \quad \mbox{and}  \quad  {\mathcal{G}} = \sqrt{\|\boldsymbol{\sigma}_0\|_{\infty}^2 + \|\boldsymbol{\sigma}_1\|_{\infty}^2G} . \label{eq:gradient_bound}
\end{align}
We will make the following four steps to obtain the desired result.
\paragraph{Step 1: Estimating A} The analysis of bounding \textbf{A} follows by introducing the proxy step-size $\asi$ and then making a decomposition as follows,
\begin{align*}
    \textbf{A} 
    &= \sum_{i=1}^d -\eta_s \cdot  \frac{\bgsi\gsi}{\bsi} + \sum_{i=1}^d \eta_s \cdot (\bgsi - \nabla f(y_s)_i) \cdot \frac{\gsi}{\bsi} \\
    &= \sum_{i=1}^d -\eta_s \cdot  \frac{\bgsi\gsi}{\asi}+\sum_{i=1}^d\eta_s \left( \frac{1}{\asi} - \frac{1}{\bsi} \right)\bgsi \gsi + \sum_{i=1}^d \eta_s \cdot (\bgsi - \nabla f(y_s)_i) \cdot \frac{\gsi}{\bsi} \\
    &=\underbrace{\sum_{i=1}^d  -\eta_s \cdot \frac{\bgsi^2}{\asi} - \sum_{i=1}^d\eta_s \cdot \frac{\bgsi \xisi}{\asi} + \sum_{i=1}^d\eta_s \left( \frac{1}{\asi} - \frac{1}{\bsi} \right)\bgsi \gsi}_{\textbf{A.1}} + \underbrace{\sum_{i=1}^d \eta_s (\bgsi - \nabla f(y_s)_i) \cdot \frac{\gsi}{\bsi}}_{\textbf{A.2}},
\end{align*}
where the last equality comes from $\xi_{s,i}=\gsi-\bgsi$.
Using the smoothness of $f$ to control the gap of $\bar{g}_{s,i}$ and $\nabla f(y_s)_i$, we have
\begin{align}
    \|\nabla f(y_s)-\bar{g}_s\| \le L\|y_s - x_s\| =  \frac{L\beta_1}{1-\beta_1}\|x_s - x_{s-1}\|. \label{eq:gradient_gap}
\end{align}
Hence, applying Young's inequality and $\eta_s \le \eta/(1-\beta_1)$ in \eqref{eq:eta_s},
\begin{align*}
    \textbf{A.2} 
    &\le \eta_s   \|\bar{g}_s - \nabla f(y_s)\| \cdot \sqrt{\sum_{i=1}^d \frac{\gsi^2}{\bsi^2}} \le \frac{1}{2L} \cdot  \|\bar{g}_s - \nabla f(y_s)\|^2 + \frac{L\eta^2}{2(1-\beta_1)^2}\cdot  \sum_{i=1}^d \frac{\gsi^2}{\bsi^2} \nonumber \\
    &\le \frac{L \beta_1^2 }{2(1-\beta_1)^2} \|x_s - x_{s-1}\|^2 + \frac{L\eta^2}{2(1-\beta_1)^2}\cdot \sum_{i=1}^d  \frac{\gsi^2}{\bsi^2}.
\end{align*}
Recalling the updated rule in Algorithm \ref{alg:Adam_scalar}, and applying \eqref{eq:hat_msi} as well as \eqref{eq:eta_s},
\begin{align}
    \textbf{A.2}&\le\frac{L \eta_{s-1}^2 \beta_1^2 }{2(1-\beta_1)^2} \sum_{i=1}^d \frac{m_{s-1,i}^2}{b_{s-1,i}^2} + \frac{L\eta^2}{2(1-\beta_1)^2}\cdot \sum_{i=1}^d \frac{\gsi^2}{\bsi^2} \nonumber \\
    &= \frac{L \eta^2 \beta_1^2 }{2(1-\beta_1)^2} \sum_{i=1}^d \frac{\hat{m}_{s-1,i}^2}{b_{s-1,i}^2} + \frac{L\eta^2}{2(1-\beta_1)^2}\cdot \sum_{i=1}^d \frac{\gsi^2}{\bsi^2}. \label{eq:A.2}
\end{align}
Then we move to estimate \textbf{A.1}. We do not directly bound \textbf{A.1}. Instead, we bound the summation of \textbf{A.1} over $s \in [t]$. First we have
\begin{align}
    \sum_{s=1}^t \textbf{A.1} =  \underbrace{- \sum_{s=1}^t \sum_{i=1}^d \eta_s \cdot \frac{\bgsi^2}{\asi}}_{\textbf{A.1.1}} \underbrace{- \sum_{s=1}^t \sum_{i=1}^d\eta_s \cdot \frac{\bgsi \xisi}{\asi}}_{\textbf{A.1.2}} + \underbrace{ \sum_{s=1}^t\sum_{i=1}^d\eta_s \left( \frac{1}{\asi} - \frac{1}{\bsi} \right)\bgsi \gsi}_{\textbf{A.1.3}}. \label{eq:(1)_decom}
\end{align}
Note that $X_{s,i} = -\frac{\eta_s \bgsi \xisi}{\asi}$ is a martingale difference sequence with respect to the filtration $\mathcal{F}_{s,i} = \sigma(\xi_{1,i},\cdots,\xi_{s,i})$ from Assumption (A2). Then setting $\omega_{s,i} = \frac{\eta_s\bgsi}{\asi}\sqrt{\sigma_{0,i}^2 + \sigma_{1,i}^2 \bgsi^2}$ and using Assumption (A3), we have
\begin{align*}
    \mE \left[ \exp\left(X_{s,i}^2 / \omega_{s,i}^2 \right) \mid \mathcal{F}_{s-1,i} \right] \le \mE \left[ \exp\left(\frac{\eta_s^2\bgsi^2 \xisi^2}{\eta_s^2\bgsi^2 (\sigma_{0,i}^2 + \sigma_{1,i}^2 \bgsi^2)} \right) \mid \mathcal{F}_{s-1,i} \right]  \le \exp(1).
\end{align*}
Noting that $\omega_{s,i}$ is $\mathcal{F}_{s-1,i}$-measurable. Hence, applying Lemma \ref{lem:Azuma} and noting that $\eta_s \le \eta/(1-\beta_1)$, we have that  for any $\lambda > 0$, with probability at least $1-\delta$, for all $t \in [T]$,
\begin{align*}
     -\sum_{s=1}^t  \eta_s \cdot \frac{\bgsi \xisi}{\asi} 
     &\le  \frac{3\lambda}{4} \sum_{s=1}^t \omega_{s,i}^2 + \frac{1}{\lambda} \log\left( \frac{T}{\delta} \right) =  \frac{3\lambda}{4}\sum_{s=1}^t  \frac{\eta_s^2\bgsi^2}{\asi^2}\left(\sigma_{0,i}^2 + \sigma_{1,i}^2 \bgsi^2 \right)+ \frac{1}{\lambda} \log\left( \frac{T}{\delta} \right) \nonumber \\
     &\le \frac{3\lambda\eta}{4(1-\beta_1)\sqrt{1-\beta_2}}\sum_{s=1}^t  \frac{\eta_s\bgsi^2}{\asi}\sqrt{\sigma_{0,i}^2 + \sigma_{1,i}^2 \bgsi^2 }+ \frac{1}{\lambda} \log\left( \frac{T}{\delta} \right),
\end{align*}
where the last inequality comes from $1/a_{s,i} \le 1/\sqrt{(1-\beta_2)(\sigma_{0,i}^2 + \sigma_{1,i}^2 \bgsi^2)}$. Hence, summing over $ i \in [d]$ and combining with \eqref{eq:gradient_bound}, we have that with probability at least $1-d\delta$,
\begin{align}
    \textbf{A.1.2}
    &\le \frac{3\lambda\eta}{4(1-\beta_1)\sqrt{1-\beta_2}}\sum_{i=1}^d\sum_{s=1}^t  \frac{\eta_s\bgsi^2}{\asi}\bar{\mathcal{G}}_s+ \frac{d}{\lambda} \log\left( \frac{T}{\delta} \right). \label{eq:1.2_middle}
\end{align}
Then setting $\lambda = (1-\beta_1)\sqrt{1-\beta_2}/\left(6\eta \mathcal{G} \right)$ and re-scaling $\delta$, we then have with probability at least $1-\delta$,
\begin{align}
    \textbf{A.1.2} \le \frac{1}{8}\sum_{i=1}^d\sum_{s=1}^t  \frac{\eta_s\bgsi^2}{\asi}{\bar{\mathcal{G}}_s \over \mathcal{G}} + \frac{6d\eta \mathcal{G}}{(1-\beta_1)\sqrt{1-\beta_2}}\log\left( \frac{dT}{\delta} \right). \label{eq:1.2}
\end{align}
Recalling Lemma \ref{lem:gap_as_bs}, then using Cauchy-Schwarz inequality, we have
\begin{align*}
    \textbf{A.1.3} &\le \sum_{i=1}^d \sum_{s=1}^t\eta_s \left| \frac{1}{\asi} - \frac{1}{\bsi} \right|\cdot  |\bgsi \gsi |\le \sum_{i=1}^d\sum_{s=1}^t \eta_s\cdot\frac{2\sqrt{1-\beta_2}}{\asi\bsi}\sqrt{\sigma_{0,i}^2 + \sigma_{1,i}^2 \bgsi^2} \cdot|\bgsi \gsi |\nonumber  \\
    &\le  \frac{1}{8}\sum_{i=1}^d\sum_{s=1}^t  \frac{\eta_s \bgsi^2}{\asi}  + 8(1-\beta_2)\sum_{i=1}^d\sum_{s=1}^t \frac{\sigma_{0,i}^2 + \sigma_{1,i}^2 \bgsi^2}{\asi} \cdot \frac{\eta_s \gsi^2}{\bsi^2} . 
\end{align*}
Applying $1/a_{s,i} \le 1/\sqrt{(1-\beta_2)(\sigma_{0,i}^2 + \sigma_{1,i}^2 \bgsi^2)}$, \eqref{eq:eta_s} and \eqref{eq:gradient_bound}, we have
\begin{align}
    \textbf{A.1.3} 
     &\le  \frac{1}{8}\sum_{i=1}^d\sum_{s=1}^t  \frac{\eta_s \bgsi^2}{\asi}  + \frac{8 \eta \sqrt{1-\beta_2}}{1-\beta_1} \sum_{i=1}^d \sum_{s=1}^t \frac{\gsi^2}{\bsi^2} \bar{\mathcal{G}}_s. \label{eq:1.3}
\end{align}
Plugging \eqref{eq:1.2} and \eqref{eq:1.3} into \eqref{eq:(1)_decom}, it then holds that with probability at least $1-\delta$,
\begin{align}
    \sum_{s=1}^t \textbf{A.1} 
    &\le - \sum_{s=1}^t \sum_{i=1}^d \left( {7\over 8} - {\bar{\mathcal{G}}_s \over 8\mathcal{G}}\right) \frac{\eta_s \bgsi^2}{\asi}+ \frac{6d\eta \mathcal{G}}{(1-\beta_1)\sqrt{1-\beta_2}}\log\left( \frac{dT}{\delta} \right) \nonumber \\
    &\quad + \frac{8 \eta \sqrt{1-\beta_2}}{1-\beta_1}\sum_{i=1}^d \sum_{s=1}^t \frac{\gsi^2}{\bsi^2} \bar{\mathcal{G}}_s. \label{eq:A.1}
 \end{align}

\paragraph{Step 2: Estimating \textbf{B}}
The analysis of bounding \textbf{B} also comes from introducing the proxy step-size $a_{s,i}$. We make a decomposition over \textbf{B},
\begin{align*}
    \textbf{B} 
    &= \underbrace{\sum_{i=1}^d \frac{\beta_1}{1-\beta_1} \cdot \left(\frac{\eta_s b_{s-1,i}}{\eta_{s-1} b_{s,i}} -1 \right)\cdot \bgsi \cdot(x_{s,i} - x_{s-1,i})}_{\textbf{B.1}} \\
    &\quad+ \underbrace{\sum_{i=1}^d \frac{\beta_1}{1-\beta_1}\cdot \left[\nabla f(y_s)_i - \bgsi\right] \cdot \left(\frac{\eta_s b_{s-1,i}}{\eta_{s-1} b_{s,i}} -1 \right)\cdot(x_{s,i} - x_{s-1,i})}_{\textbf{B.2}}.
\end{align*}
First applying the updated rule,
\begin{align*}
    \textbf{B.1} &\le \sum_{i=1}^d \frac{\beta_1}{1-\beta_1} \cdot \left|\left(\frac{\eta_s b_{s-1,i}}{\eta_{s-1} b_{s,i}} -1 \right)\cdot \bgsi \cdot \frac{\eta_{s-1}m_{s-1,i}}{b_{s-1,i}} \right| \\
    &= \sum_{i=1}^d \frac{\beta_1}{1-\beta_1} \cdot \left|\left(\frac{\eta_s }{ b_{s,i}} -\frac{\eta_{s-1} }{ b_{s-1,i}}\right)\cdot \bgsi \cdot m_{s-1,i} \right| \nonumber \\
    &\le \sum_{i=1}^d \frac{\beta_1}{1-\beta_1} \cdot \left|\left(\frac{\eta_s }{ b_{s,i}} -\frac{\eta_{s} }{ b_{s-1,i}}\right)\cdot \bgsi \cdot m_{s-1,i} \right| + \sum_{i=1}^d \frac{\beta_1}{1-\beta_1} \cdot \left|\left(\eta_{s-1}-\eta_s\right) \cdot \bgsi \cdot \frac{m_{s-1,i}}{b_{s-1,i}} \right|. \nonumber
\end{align*}
It's hard to well estimate the gap of $|1/\bsi-1/b_{s-1,i}|$ due to the exponential moving average. Hence, we introduce $a_{s,i}$ and obtain the following decomposition 
\begin{align}
    \textbf{B.1} 
    &\le \underbrace{\sum_{i=1}^d \frac{\beta_1}{1-\beta_1} \cdot \left|\frac{1}{b_{s,i}} -\frac{1}{a_{s,i}}\right| \cdot \eta_{s}\cdot |\bgsi \cdot m_{s-1,i}|}_{\textbf{B.1.1}} + \underbrace{\sum_{i=1}^d \frac{\beta_1}{1-\beta_1} \cdot \left|\frac{1}{a_{s,i}} -\frac{1}{b_{s-1,i}}\right|  \cdot \eta_{s}\cdot |\bgsi \cdot m_{s-1,i}|}_{\textbf{B.1.2}} \nonumber \\
    &+ \underbrace{\sum_{i=1}^d \frac{\beta_1}{1-\beta_1} \cdot \left|\left(\eta_{s-1}-\eta_s\right) \cdot \bgsi \cdot \frac{m_{s-1,i}}{b_{s-1,i}} \right|}_{\textbf{B.1.3}}. \label{eq:B.1_decom}
\end{align}
Then applying Lemma \ref{lem:gap_as_bs} and Young's inequality, and using $1/\asi \le 1/\sqrt{(1-\beta_2)(\sigma_{0,i}^2 + \sigma_{1,i}^2 \bgsi^2)}$,
\begin{align}
    \textbf{B.1.1} 
    &\le \sum_{i=1}^d\frac{\beta_1}{1-\beta_1}\cdot \frac{2\eta_{s} \sqrt{1-\beta_2}}{\asi \bsi} \cdot \left(\sqrt{\sigma_{0,i}^2 + \sigma_{1,i}^2 \bgsi^2}\right) |\bgsi \cdot m_{s-1,i}| \nonumber \\
    &\le \sum_{i=1}^d\frac{\eta_s}{4} \cdot \frac{\bgsi^2}{\asi} +  \frac{4 \eta_s \beta_1^2 (1-\beta_2)}{(1-\beta_1)^2}\cdot \sum_{i=1}^d\frac{\sigma_{0,i}^2 + \sigma_{1,i}^2 \bgsi^2}{\asi} \cdot \frac{m_{s-1,i}^2}{\bsi^2}\nonumber  \\
    &\le \sum_{i=1}^d\frac{\eta_s}{4} \cdot \frac{\bgsi^2}{\asi} + \frac{4\eta_s \beta_1^2 \sqrt{1-\beta_2}}{(1-\beta_1)^2}\cdot  \sum_{i=1}^d \left(\sqrt{\sigma_{0,i}^2 + \sigma_{1,i}^2 \bgsi^2} \right)\frac{m_{s-1,i}^2}{\bsi^2}.\label{eq:B.1.1_middle}
\end{align}
Further applying \eqref{eq:eta_s} and \eqref{eq:gradient_bound}, we also have 
\begin{align}
    \textbf{B.1.1} 
    &\le \sum_{i=1}^d\frac{\eta_s}{4} \cdot \frac{\bgsi^2}{\asi} + \frac{4\eta_s \beta_1^2 \sqrt{1-\beta_2}}{(1-\beta_1)^2} \sum_{i=1}^d\frac{m_{s-1,i}^2}{\bsi^2} \bar{\mathcal{G}}_s \nonumber \\
    &\le \sum_{i=1}^d\frac{\eta_s}{4} \cdot \frac{\bgsi^2}{\asi} + \left(c_1 + c_2\sqrt{G_t}\right)\cdot\sum_{i=1}^d\frac{m_{s-1,i}^2}{\bsi^2}, \label{eq:B.1.1}
\end{align}
where we denote the constants $c_1,c_2$ as
\begin{align}\label{eq:define_c_1}
    c_1 = \frac{4\|\boldsymbol{\sigma}_0\|_{\infty} \eta \beta_1^2 \sqrt{1-\beta_2}}{(1-\beta_1)^3} , \quad c_2 = \frac{4\|\boldsymbol{\sigma}_1\|_{\infty} \eta \beta_1^2 \sqrt{1-\beta_2}}{(1-\beta_1)^3}.
\end{align}
Similarly applying Lemma \ref{lem:gap_as_bs} and Young's inequality, we have
\begin{align*}
    \textbf{B.1.2} 
    &\le \sum_{i=1}^d\frac{\eta_{s}\beta_1}{1-\beta_1}\frac{\sqrt{1-\beta_2}}{b_{s-1,i}\asi}\cdot \max\left\{ \sqrt{\sigma_{0,i}^2 + (1+\sigma_{1,i}^2)\bgsi^2} + \ep,  \sqrt{v_{s-1,i}} + \ep \right\}\cdot |\bgsi \cdot m_{s-1,i}| \\
    &\le \sum_{i=1}^d \frac{\eta_s}{4}\cdot \frac{\bgsi^2}{\asi} + \frac{\eta_s \beta_1^2(1-\beta_2)}{(1-\beta_1)^2} \\
    &\quad\cdot \sum_{i=1}^d \frac{1}{\asi}\max\left\{ \left(\sqrt{\sigma_{0,i}^2 + (1+\sigma_{1,i}^2)\bgsi^2} + \ep\right)^2,  \left(\sqrt{v_{s-1,i}} + \ep\right)^2 \right\} \cdot \frac{m_{s-1,i}^2}{b_{s-1,i}^2}.
\end{align*}
Noting that since $1-\beta_2 < 1-\beta_2^s$ and $\beta_2 \in (0,1)$ and recalling \eqref{eq:prox_as}, we have
\begin{align*}
    \frac{1}{\asi} 
    &\le \frac{1}{\sqrt{(1-\beta_2)(\sigma_{0,i}^2 + (1+\sigma_{1,i}^2)\bgsi^2)} + \ep\sqrt{1-\beta_2^s}} \le 
    \frac{1}{\sqrt{1-\beta_2} \left(\sqrt{\sigma_{0,i}^2 + (1+\sigma_{1,i}^2)\bgsi^2}+ \ep \right) } \\
    &\le  \frac{1}{\sqrt{\beta_2(1-\beta_2)} \left(\sqrt{\sigma_{0,i}^2 + (1+\sigma_{1,i}^2)\bgsi^2}+ \ep \right)}. 
\end{align*}
We also have
\begin{align*}
    \frac{1}{\asi} \le \frac{1}{\sqrt{\beta_2 v_{s-1,i}} + \ep\sqrt{1-\beta_2^s}} \le \frac{1}{\sqrt{\beta_2(1-\beta_2)} \left(\sqrt{v_{s-1,i}}+ \ep \right) }.
\end{align*}
Hence, we further obtain 
\begin{align}
    \textbf{B.1.2} 
    &\le \sum_{i=1}^d \frac{\eta_s}{4}\cdot \frac{\bgsi^2}{\asi} + \frac{\eta_s \beta_1^2\sqrt{1-\beta_2}}{(1-\beta_1)^2\sqrt{\beta_2}} \nonumber\\
    &\quad\cdot \sum_{i=1}^d\max\left\{ \sqrt{\sigma_{0,i}^2 + (1+\sigma_{1,i}^2)\bgsi^2} + \ep,  \sqrt{v_{s-1,i}} + \ep \right\} \cdot \frac{m_{s-1,i}^2}{b_{s-1,i}^2}. \label{eq:B.1.2_middle}
\end{align}
Applying the definition of $v_{s-1,i}$ from Algorithm \ref{alg:Adam_scalar} and Assumption (A3), 
\begin{align*}
    \sqrt{v_{s-1,i}} 
    &= \sqrt{(1-\beta_2)\sum_{j=1}^{s-1}\beta_2^{s-1-j}g_{j,i}^2} \le \sqrt{2(1-\beta_2)\sum_{j=1}^{s-1}\beta_2^{s-1-j}\left(\bar{g}_{j,i}^2 + \xi_{j,i}^2\right)} \\
    &\le \sqrt{2(1-\beta_2)\sum_{j=1}^{s-1}\beta_2^{s-1-j}\left[  \sigma_{0,i}^2 + (1+\sigma_{1,i}^2)\bar{g}_{j,i}^2\right]} 
   \le \sqrt{2\left[\sigma_{0,i}^2 + (1+\sigma_{1,i}^2) \left(\max_{j \in [s]} \bar{g}_{j,i}^2 \right) \right]}.
\end{align*}
As a consequence of \eqref{eq:gradient_bound}, we have $\max_{j \in [s]} \bar{g}_{j,i}^2 \le G_t, \forall s \in [t]$. Hence, we further have 
\begin{align}
    \max\left\{ \sqrt{\sigma_{0,i}^2 + (1+\sigma_{1,i}^2)\bgsi^2} + \ep,  \sqrt{v_{s-1,i}} + \ep \right\} 
    &\le \sqrt{2\left[\|\boldsymbol{\sigma}_0\|_{\infty}^2 + (1+\|\boldsymbol{\sigma}_1\|_{\infty}^2)G_t \right]} + \ep. \label{eq:max_1_middle}
\end{align}
Hence, plugging \eqref{eq:max_1_middle} into \eqref{eq:B.1.2_middle}, and applying \eqref{eq:eta_s} and $\ep = \ep_0\sqrt{1-\beta_2 }\le \ep_0$, 
\begin{align}
     \textbf{B.1.2} 
     &\le  \sum_{i=1}^d \frac{\eta_s}{4}\cdot \frac{\bgsi^2}{\asi} + \frac{\eta \beta_1^2\sqrt{1-\beta_2}}{(1-\beta_1)^3\sqrt{\beta_2}}\cdot\left( \ep + \sqrt{2\|\boldsymbol{\sigma}_0\|_{\infty}^2 +  2(1+\|\boldsymbol{\sigma}_1\|_{\infty}^2)G_t}  \right) \sum_{i=1}^d  \frac{m_{s-1,i}^2}{b_{s-1,i}^2}  \nonumber \\
     &\le \sum_{i=1}^d \frac{\eta_s}{4}\cdot \frac{\bgsi^2}{\asi} +  \left(d_1 + d_2\sqrt{G_t}\right) \cdot \sum_{i=1}^d\frac{m_{s-1,i}^2}{b_{s-1,i}^2} . \label{eq:B.1.2}
\end{align}
where we denote constants $d_1, d_2$ as
\begin{align}\label{eq:define_d_1}
    d_1 = \frac{\eta\beta_1^2\sqrt{1-\beta_2}}{(1-\beta_1)^3\sqrt{\beta_2}}\cdot\left( \ep_0 + \sqrt{2}\|\boldsymbol{\sigma}_0\|_{\infty}  \right), \quad d_2=   \frac{\eta\beta_1^2\sqrt{2(1-\beta_2)(1+\|\boldsymbol{\sigma}_1\|_{\infty}^2)}}{(1-\beta_1)^3\sqrt{\beta_2}}.
\end{align}
We then plug \eqref{eq:B.1.1} and \eqref{eq:B.1.2} into \eqref{eq:B.1_decom},
\begin{align}
    \textbf{B.1} 
    &\le \sum_{i=1}^d \frac{\eta_s}{2}\cdot \frac{\bgsi^2}{\asi} + \left(c_1 + c_2\sqrt{G_t} \right) \sum_{i=1}^d \frac{m_{s-1,i}^2}{b_{s,i}^2} + \left(d_1 + d_2\sqrt{G_t}\right) \sum_{i=1}^d\frac{m_{s-1,i}^2}{b_{s-1,i}^2} + \textbf{B.1.3}. \label{eq:B.1}
\end{align}
Then we move to bound the summation of \textbf{B.1.3} over $s \in \{2,\cdots, t \}$ since $m_{0}=0$. Recalling \eqref{eq:eta_s}, we have the following estimation,
\begin{align}
    \textbf{B.1.3} 
    &\le \underbrace{\sum_{i=1}^d \frac{\eta \beta_1\sqrt{1-\beta_2^s}}{1-\beta_1} \cdot \left|\left(\frac{1}{1-\beta_1^{s-1}}-\frac{1}{1-\beta_1^{s}}\right) \cdot \bgsi \cdot \frac{m_{s-1,i}}{b_{s-1,i}} \right|}_{\textbf{B.1.3.1}} \nonumber\\
    &\quad+ \underbrace{\sum_{i=1}^d \frac{\eta \beta_1}{(1-\beta_1)(1-\beta_1^{s-1})} \cdot \left|\left(\sqrt{1-\beta_2^{s-1}} - \sqrt{1-\beta_2^{s}}\right) \cdot \bgsi  \cdot \frac{m_{s-1,i}}{b_{s-1,i}} \right|}_{\textbf{B.1.3.2}}. \label{eq:B.1.3_decomp}
\end{align}
We first combine with Lemma \ref{lem:OrderT_bound_gradient} and \eqref{eq:gradient_bound}, 
\begin{align}
    \sqrt{1-\beta_2^s}\left|\frac{\bgsi m_{s-1,i}}{b_{s-1,i}} \right|\le \sqrt{1-\beta_2^s}\|\bar{g}_s\|_{\infty} \cdot \left|\frac{m_{s-1,i}}{b_{s-1,i}} \right| \le \sqrt{G_t}\cdot \sqrt{\frac{(1-\beta_1)(1-\beta_1^{s-1})}{(1-\beta_2)(1-\beta_1/\beta_2)}}. \label{eq:bound_1}
\end{align}
Noting that $1/(1-\beta_1^{s-1}) \ge 1/(1-\beta_1^{s})$ for any $s \ge 2$, thus we have for any $s \ge 2$,
\begin{align}
    \textbf{B.1.3.1} \le \frac{d\eta \beta_1}{1-\beta_1}\sqrt{G_t}\cdot \sqrt{\frac{(1-\beta_1)(1-\beta_1^{s-1})}{(1-\beta_2)(1-\beta_1/\beta_2)}}\cdot\left(\frac{1}{1-\beta_1^{s-1}}-\frac{1}{1-\beta_1^{s}}\right) . \label{eq:B.1.3.1}
\end{align}
Hence, summing over $s \in [t]$ on \textbf{B.1.3.1} and noting that \textbf{B.1.3.1} vanishes when $s = 1$, 
\begin{align}
    \sum_{s=1}^t \textbf{B.1.3.1} 
    &\le \frac{d\eta\beta_1}{1-\beta_1}\sqrt{G_t}\cdot \sqrt{\frac{(1-\beta_1)(1-\beta_1^{s-1})}{(1-\beta_2)(1-\beta_1/\beta_2)}}\sum_{s=2}^t \left(\frac{1}{1-\beta_1^{s-1}}-\frac{1}{1-\beta_1^{s}}\right) \nonumber\\
    &\le d\eta \beta_1 \cdot \sqrt{\frac{1}{(1-\beta_1)(1-\beta_2)(1-\beta_1/\beta_2)}}\cdot \left( \frac{1}{1-\beta_1} - \frac{1}{1-\beta_1^{t}} \right)\sqrt{G_t} \nonumber \\
    &\le \frac{d\eta \beta_1 }{\sqrt{(1-\beta_1)^3(1-\beta_2)(1-\beta_1/\beta_2)}}\sqrt{G_t}:=e\sqrt{G_t}.
    \label{eq:sum_B.1.3.1}
\end{align}
We also combine with Lemma \ref{lem:OrderT_bound_gradient}, \eqref{eq:eta_s} and \eqref{eq:gradient_bound},
\begin{align*}
    \frac{1}{1-\beta_1^{s-1}}\left|\frac{\bgsi m_{s-1,i}}{b_{s-1,i}} \right| 
     \le \sqrt{G_t}\cdot \sqrt{\frac{1}{(1-\beta_2)(1-\beta_1/\beta_2)}}.
\end{align*}
Hence, summing over $s \in [t]$ on \textbf{B.1.3.2} and noting that \textbf{B.1.3.2} vanishes when $s = 1$, 
\begin{align}
    \sum_{s=1}^t \textbf{B.1.3.2} 
    &\le \frac{d\eta\beta_1}{1-\beta_1}\sqrt{G_t}\cdot \sqrt{\frac{1}{(1-\beta_2)(1-\beta_1/\beta_2)}} \sum_{s=2}^t \left(\sqrt{1-\beta_2^{s}} - \sqrt{1-\beta_2^{s-1}}\right) \nonumber \\
    &\le  \frac{d\eta \beta_1 }{(1-\beta_1)\sqrt{(1-\beta_2)(1-\beta_1/\beta_2)}}\sqrt{G_t}:=f\sqrt{G_t}. \label{eq:sum_B.1.3.2}
\end{align}
Then we move to upper bound $\textbf{B.2}$. Note that when $s = 1$, \textbf{B.2} vanishes since $x_0 = x_1$ from the beginning. Now we consider $s \ge 2$. When $\frac{\eta_s b_{s-1,i}}{\eta_{s-1} b_{s,i}} < 1$, we have 
\begin{align*}
    \left|\frac{\eta_s b_{s-1,i}}{\eta_{s-1} \bsi}-1\right| = 1 - \frac{\eta_s b_{s-1}}{\eta_{s-1} b_s} < 1.
\end{align*}
Let $p = \beta_2^{s-1}$, since $0 < 1-\beta_1^{s-1} < 1-\beta_1^s$, then we have
\begin{align*}
    \frac{\eta_s}{\eta_{s-1}} = \sqrt{\frac{1-\beta_2^s}{1-\beta_2^{s-1}}}\cdot \frac{1-\beta_1^{s-1}}{1-\beta_1^s} \le  \sqrt{1+\frac{\beta_2^{s-1}(1-\beta_2)}{1-\beta_2^{s-1}}} = \sqrt{1+(1-\beta_2)\cdot \frac{p}{1-p}}.
\end{align*}
Since $h(p)=p/(1-p)$ is increasing as $p$ grows, we then obtain that $p$ takes the maximum value when $s = 2$. Hence, it holds that
\begin{align}
     \frac{\eta_s}{\eta_{s-1}} \le \sqrt{1+(1-\beta_2)\cdot \frac{\beta_2}{1-\beta_2}} = \sqrt{1+\beta_2}. \label{eq:Sigma_max_1}
\end{align}
Then, since $\ep_{s-1} \le \ep_s$, we further have 
\begin{align}
    \frac{b_{s-1,i}}{\bsi} = \frac{\ep_{s-1} + \sqrt{v_{s-1,i}} }{\ep_s + \sqrt{\beta_2v_{s-1,i}+(1-\beta_2)\gsi^2 }} \le \frac{\ep_{s} + \sqrt{v_{s-1,i}} }{\ep_s + \sqrt{\beta_2v_{s-1,i}}} \le \frac{1}{\sqrt{\beta_2}}. \label{eq:Sigma_max_2}
\end{align}
Combining with \eqref{eq:Sigma_max_1} and \eqref{eq:Sigma_max_2}, we have when $ \frac{\eta_s b_{s-1,i}}{\eta_{s-1}\bsi} \ge 1$
\begin{align*}
     \left|\frac{\eta_s b_{s-1,i}}{\eta_{s-1} \bsi}-1\right| = \frac{\eta_s b_{s-1,i}}{\eta_{s-1} \bsi}-1 \le \sqrt{\frac{1+\beta_2}{\beta_2}} -1 .
\end{align*}
Hence, we have for any $i \in [d]$ and any $s \ge 2$,
\begin{align}
    \left|\frac{\eta_s b_{s-1,i}}{\eta_{s-1} b_{s,i}}-1\right| \le \Sigma_{\max}:=\max\left\{1, \sqrt{\frac{1+\beta_2}{\beta_2}} -1  \right\}. \label{eq:Sigma_max}
\end{align}
Then, we apply the H\"{o}lder's inequality and recalling \eqref{eq:gradient_gap},
\begin{align*}
    \textbf{B.2 } 
    &\le \frac{\beta_1\Sigma_{\max}}{1-\beta_1}\cdot \sum_{i=1}^d  \left|\nabla f(y_s)_i - \bgsi\right| \cdot | x_{s,i} - x_{s-1,i}| \\
    &\leq \frac{\beta_1\Sigma_{\max}}{1-\beta_1}\cdot 
    \|\nabla f(y_s) - \bar{g}_s\| \cdot \|x_s - x_{s-1}\| \\
    &\le \frac{L\beta_1^2\Sigma_{\max}^2}{(1-\beta_1)^2}\cdot \|x_s - x_{s-1}\|^2,
\end{align*}
where the last inequality comes from $\Sigma_{\max} \ge 1$. Further using the updated rule and \eqref{eq:eta_s}, and recalling \eqref{eq:hat_msi}, we further obtain that
\begin{align}
   \textbf{B.2} 
   \le \frac{L\beta_1^2\Sigma_{\max}^2}{(1-\beta_1)^2}\cdot  \eta_{s-1}^2 \sum_{i=1}^d \frac{m_{s-1,i}^2}{b_{s-1,i}^2} \le \frac{L\eta^2\beta_1^2\Sigma_{\max}^2}{(1-\beta_1)^2}\cdot \sum_{i=1}^d \frac{\hat{m}_{s-1,i}^2}{b_{s-1,i}^2}. \label{eq:B.2}
\end{align}
\paragraph{Step 3: Estimating \textbf{C}}
Applying Young's inequality, \eqref{eq:Sigma_max} and the updated rule, and combining with \eqref{eq:hat_msi} and \eqref{eq:eta_s}, 
\begin{align}
    \textbf{C} 
    &\le L\eta_s^2 \sum_{i=1}^d \frac{\gsi^2}{\bsi^2} + \frac{L\beta_1^2\Sigma_{\max}^2}{(1-\beta_1)^2} \|x_s - x_{s-1}\|^2 \nonumber \\
    &\le L\eta_s^2 \sum_{i=1}^d \frac{\gsi^2}{\bsi^2} + \frac{L\eta_{s-1}^2\beta_1^2\Sigma_{\max}^2}{(1-\beta_1)^2} \sum_{i=1}^d  \frac{m_{s-1,i}^2}{b_{s-1,i}^2} \nonumber \\
    &\le \frac{L\eta^2}{(1-\beta_1)^2} \sum_{i=1}^d \frac{\gsi^2}{\bsi^2} + \frac{L\eta^2\beta_1^2\Sigma_{\max}^2}{(1-\beta_1)^2} \sum_{i=1}^d  \frac{\hat{m}_{s-1,i}^2}{b_{s-1,i}^2}\label{eq:C}.
\end{align}
\paragraph{Step 4: Putting together} Now we have separately estimate all terms in \eqref{eq:A+B+C}. Note that combining with Lemma \ref{lem:gradient_delta_s} and Lemma \ref{lem:gradient_xs_ys}, we have
\begin{align}
    \|\bar{g}_{t+1}\|^2 \le 2\|\nabla f(y_{t+1})\|^2 + 2M^2 \le 4L(f(y_{t+1}) - f^*) + 2M^2. \label{eq:LHS}
\end{align}
We will now prove the result by induction. First from Lemma \ref{lem:gradient_delta_s}, it's easy to obtain that $G_1 = \|\bar{g}_1\|^2 \le 2L(f(x_1)-f^*) \le G$ since $x_1=y_1$ from \eqref{eq:define_y_s}. Then we assume that
\begin{align}\label{eq:assum_boundgradient}
	G_s \leq  G,\quad \forall s \in [t].
\end{align}
Then we obtain $\bar{\mathcal{G}}_s \leq \mathcal{G}$ from \eqref{eq:gradient_bound}. Hence, combining with \eqref{eq:A.1}, 
\begin{align}
    \sum_{s=1}^t \textbf{A.1} 
    &\le -\frac{3}{4} \sum_{s=1}^t \sum_{i=1}^d  \frac{\eta_s \bgsi^2}{\asi}+ \frac{6d\eta \mathcal{G}}{(1-\beta_1)\sqrt{1-\beta_2}}\log\left( \frac{dT}{\delta} \right)  + \frac{8 \eta\mathcal{G} \sqrt{1-\beta_2}}{1-\beta_1}\sum_{i=1}^d \sum_{s=1}^t \frac{\gsi^2}{\bsi^2} . \label{eq:A.1_new}
 \end{align}
Then sum up over $s\in [t]$ in \eqref{eq:A+B+C} and then plug \eqref{eq:A.2}, \eqref{eq:B.1}, \eqref{eq:sum_B.1.3.1}, \eqref{eq:sum_B.1.3.2}, \eqref{eq:B.2}, \eqref{eq:C} and \eqref{eq:A.1_new} into it. With the inducted assumption \eqref{eq:assum_boundgradient} and \eqref{eq:gradient_bound}, we have 
\begin{align}
    f(y_{t+1}) &\le f(y_1) - \sum_{s=1}^t \sum_{i=1}^d \frac{\eta_s}{4} \cdot \frac{\bgsi^2}{\asi} + (e+f)\sqrt{G} +\frac{6d\eta\sqrt{\|\boldsymbol{\sigma}_0\|_{\infty}^2 + \|\boldsymbol{\sigma}_1\|_{\infty}^2G}}{(1-\beta_1)\sqrt{1-\beta_2}}\log\left( \frac{dT}{\delta} \right)\nonumber \\
    &\quad+ \left( \frac{8 \eta \sqrt{1-\beta_2}}{1-\beta_1}\cdot \sqrt{\|\boldsymbol{\sigma}_0\|_{\infty}^2 +  \|\boldsymbol{\sigma}_1\|_{\infty}^2 G }+\frac{3L\eta^2}{2(1-\beta_1)^2}  \right) \sum_{i=1}^d \sum_{s=1}^t \frac{\gsi^2}{\bsi^2}
    \nonumber \\
    &\quad 
    + \left(c_1 + c_2 \sqrt{G} \right)\sum_{s=1}^{t-1} \sum_{i=1}^d  \frac{m_{s,i}^2}{b_{s+1,i}^2} +\left( d_1 + d_2\sqrt{G} \right)\sum_{s=1}^{t-1} \sum_{i=1}^d  \frac{m_{s,i}^2}{b_{s,i}^2} \nonumber \\
    &\quad+\left( \frac{L\eta^2\beta_1^2}{2(1-\beta_1)^2} + \frac{2L\eta^2\beta_1^2\Sigma_{\max}^2}{(1-\beta_1)^2} \right)\sum_{s=1}^{t-1} \sum_{i=1}^d  \frac{\hat{m}_{s,i}^2}{b_{s,i}^2}. \label{eq:(1)+(2)+(3)}
\end{align}
We then combine with Lemma \ref{lem:sum_1} and $\sqrt{a + b G} \leq \sqrt{a} + \sqrt{b G} \leq \sqrt{a} + {cG\over 2} + {G \over 2c}$
for any $a,b,c>0$, and apply $\eta = C_0\sqrt{1-\beta_2}$, 
\begin{align}
    f(y_{t+1}) &\le \frac{G}{8L} + f(y_1)- \sum_{s=1}^t \sum_{i=1}^d \frac{\eta_s}{4} \cdot \frac{\bgsi^2}{\asi}  +16L(e^2+f^2)  \nonumber \\
    &\quad+ \frac{6dC_0\|\boldsymbol{\sigma}_0\|_{\infty} }{1-\beta_1}\log\left(\frac{dT}{\delta} \right)+ \frac{288Ld^2C_0^2\|\boldsymbol{\sigma}_1\|_{\infty}^2 }{(1-\beta_1)^2}\log^2\left(\frac{dT}{\delta} \right) \nonumber \\
    &\quad +\left(\frac{3LC_0^2}{2(1-\beta_1)^2}   + \frac{8\|\boldsymbol{\sigma}_0\|_{\infty}C_0}{1-\beta_1}\right)d\log\left( \frac{\mathcal{F}(T)}{\beta_2^T} \right)+  \frac{512L d^2C_0^2\|\boldsymbol{\sigma}_1\|_{\infty}^2}{(1-\beta_1)^2}\log^2\left( \frac{\mathcal{F}(T)}{\beta_2^T} \right)
     \nonumber \\
    &\quad+ \frac{c_1 (1-\beta_1)d}{\beta_2(1-\beta_2)(1-\beta_1/\beta_2)}\log\left( \frac{\mathcal{F}(T)}{\beta_2^T} \right)+ \frac{16Lc_2^2(1-\beta_1)^2d^2}{\beta_2^2(1-\beta_2)^2(1-\beta_1/\beta_2)^2}\log^2\left( \frac{\mathcal{F}(T)}{\beta_2^T} \right)  \nonumber \\
    &\quad +\left( \frac{LC_0^2\beta_1^2}{2(1-\beta_1)^2} + \frac{2LC_0^2\beta_1^2\Sigma^2_{\max}}{(1-\beta_1)^2}  \right)\frac{d }{1-\beta_1/\beta_2}\log\left( \frac{\mathcal{F}(T)}{\beta_2^T} \right)\nonumber\\
    &\quad+\frac{d_1(1-\beta_1)d}{(1-\beta_2)(1-\beta_1/\beta_2)}\log\left( \frac{\mathcal{F}(T)}{\beta_2^T} \right)  + \frac{16Ld_2^2(1-\beta_1)^2d^2}{(1-\beta_2)^2(1-\beta_1/\beta_2)^2}\log^2\left( \frac{\mathcal{F}(T)}{\beta_2^T} \right) . \label{eq:delta_s_final_1}
\end{align}
Then with both sides of  \eqref{eq:delta_s_final_1} subtracting $f^*$ and then combining with \eqref{eq:LHS}, we could obtain an upper bound for $\|\bar{g}_{t+1}\|^2$. Further combining with the definition of $c_1,c_2,d_1,d_2,e,f, M$ in \eqref{eq:define_c_1}, \eqref{eq:define_d_1}, \eqref{eq:sum_B.1.3.1}, \eqref{eq:sum_B.1.3.2} and Lemma \ref{lem:gradient_xs_ys}, we obtain the desired bound \eqref{eq:delta_s_final_2} holds with probability at least $1-\delta$,
where $G$ is given by the following expression that is independent from $t$. Note that we apply $\beta_1,\beta_2 \in [0,1)$ and $x_1=y_1$ from \eqref{eq:define_y_s} to simplify the expression with details omitted.
\begin{align}\label{eq:G_define}
    G &:= 8L(f(x_1)-f^*)  +\frac{256L^2C_0^2 \beta_1^2d^2}{(1-\beta_1)^3(1-\beta_1/\beta_2)}
    +\frac{4L^2C_0^2\beta_1^2d}{(1-\beta_1)^2(1-\beta_1/\beta_2)}
    +\frac{48LC_0\|\boldsymbol{\sigma}_0\|_{\infty} }{1-\beta_1}d\log\left( \frac{dT}{\delta} \right)
    \nonumber \\
    &\quad+ \frac{8LC_0\beta_1^2}{(1-\beta_1)^2(1-\beta_1/\beta_2)}\left[\frac{4\|\boldsymbol{\sigma}_0\|_{\infty} }{\beta_2} +  \frac{\ep_0 + \sqrt{2}\|\boldsymbol{\sigma}_0\|_{\infty}}{\sqrt{\beta_2}}+\frac{LC_0}{2}(1+4\Sigma_{\max}^2)\right] d\log\left(\frac{\mathcal{F}(T)}{\beta_2^T}\right)  \nonumber \\
    &\quad+ 8\left(\frac{3L^2C_0^2}{(1-\beta_1)^2}  + \frac{8L\|\boldsymbol{\sigma}_0\|_{\infty} C_0}{1-\beta_1} \right)d\log\left(\frac{\mathcal{F}(T)}{\beta_2^T}\right)+ \frac{2304L^2C_0^2\|\boldsymbol{\sigma}_1\|_{\infty}^2 }{(1-\beta_1)^2}d^2\log^2\left( \frac{dT}{\delta} \right) \nonumber \\
    &\quad+256\left[\frac{16L^2C_0^2\|\boldsymbol{\sigma}_1\|_{\infty}^2 }{(1-\beta_1)^2}+ \frac{9L^2C_0^2\beta_1^4 (1+\|\boldsymbol{\sigma}_1\|_{\infty}^2)}{\beta_2^2(1-\beta_1)^4(1-\beta_1/\beta_2)^2}\right]d^2\log^2\left(\frac{\mathcal{F}(T)}{\beta_2^T}\right)
    .
\end{align}
Finally from \eqref{eq:delta_s_final_2} we obtain that $\|\nabla f(x_{t+1})\|^2 \le G$. Combining with \eqref{eq:assum_boundgradient}, the induction is completed and we obtain $\|\nabla f(x_t)\|^2 \le G, \forall t \in [T+1]$.
\hfill \BlackBox

\subsection{Proof of Theorem \ref{thm:Adam_scalar_1}}
    Let us set $t = T$ in \eqref{eq:delta_s_final_2}. We then re-range the order and have that with probability at least $1-\delta$,
    \begin{align}
        L \sum_{s=1}^T \sum_{i=1}^d \eta_s \cdot \frac{\bgsi^2}{\asi}  &\le G - \|\bar{g}_{T+1}\|^2 \le G. \label{eq:final_middle_4}
    \end{align}
    Recalling the definition of $\asi=\ep_s + \sqrt{\tilde{v}_{s,i}}$ in \eqref{eq:prox_as}, we denote $\tilde{\mathbf{v}}_{s} = (\tilde{v}_{s,i})_i \in \left(\mR^+ \right)^d$, then
    \begin{align}
        G \ge L \sum_{s=1}^T \sum_{i=1}^d \eta_s \cdot \frac{\bgsi^2}{\ep_s + \sqrt{\tilde{v}_{s,i} }}
        \ge  L \sum_{s=1}^T \sum_{i=1}^d \eta_s \cdot \frac{\bgsi^2}{\ep_s + \sqrt{\|\tilde{\mathbf{v}}_{s}\|_1} } = L \sum_{s=1}^T \eta_s \cdot \frac{\|\bar{g}_s\|^2}{\ep_s + \sqrt{\|\tilde{\mathbf{v}}_{s}\|_1}} . \label{eq:final_middle_5}
    \end{align}
    Then applying the definition of $\tilde{v}_{s,i}$, the basic inequality and Assumption (A3), we have
    \begin{align}
        \sqrt{\|\tilde{\mathbf{v}}_{s}\|_1 } 
        &= \sqrt{\sum_{i=1}^d |\tilde{v}_{s,i}| } =  \sqrt{\sum_{i=1}^d  (1-\beta_2) \left[\sum_{j=1}^{s-1} \beta_2^{s - j} g_{j,i}^2+ (\sigma_{0,i}^2 + (1+\sigma_{1,i}^2)\bar{g}_{s,i}^2)\right]  } \nonumber \\
        &\le \sqrt{\sum_{i=1}^d 2(1-\beta_2) \left[\sum_{j=1}^{s-1} \beta_2^{s - j} (\sigma_{0,i}^2 + (1+\sigma_{1,i}^2)\bar{g}_{j,i}^2)+ (\sigma_{0,i}^2 + (1+\sigma_{1,i}^2)\bar{g}_{s,i}^2)\right] } \nonumber \\
        &\le \sqrt{2(1-\beta_2)\sum_{i=1}^d  \sum_{j=1}^{s}\beta_2^{s - j} (\sigma_{0,i}^2 + (1+\|\boldsymbol{\sigma}_1\|_{\infty}^2) \bar{g}_{j,i}^2) } \nonumber \\
        &= \sqrt{2(1-\beta_2)\left(\|\boldsymbol{\sigma}_0\|^2 \sum_{j=1}^{s} \beta_2^{s - j} + (1+\|\boldsymbol{\sigma}_1\|_{\infty}^2)  \sum_{j=1}^{s} \beta_2^{s - j} \|\bar{g}_j\|^2 \right) }. \label{eq:bound_vs}
    \end{align}
    Recalling the upper bound $\|\bar{g}_j\|^2 \le G, \forall j \in [T+1]$ in Proposition \ref{pro:delta_s_1}. Hence, we have with probability at least $1-\delta$,
    \begin{align*}
        \sqrt{\|\tilde{\mathbf{v}}_{s}\|_1 } 
        &\le \sqrt{2(1-\beta_2)\left(\|\boldsymbol{\sigma}_0\|^2 \sum_{j=1}^{s} \beta_2^{s - j} + (1+\|\boldsymbol{\sigma}_1\|_{\infty}^2)  G\sum_{j=1}^{s} \beta_2^{s - j} \right) }  \\
        &\le \sqrt{2(1-\beta_2^s)\left(\|\boldsymbol{\sigma}_0\|^2 + (1+\|\boldsymbol{\sigma}_1\|_{\infty}^2)G \right)   }.
    \end{align*}
    From the setting $\eta_s = C_0 \sqrt{(1-\beta_2^s)(1-\beta_2)}/(1-\beta_1^s)$ and $\ep_s = \ep_0\sqrt{(1-\beta_2^s)(1-\beta_2)}$ in \eqref{eq:parameter_setting}, we then have for $s \in [T]$,
    \begin{align*}
        \frac{\ep_s +\sqrt{\|\tilde{\mathbf{v}}_{s}\|_1 }  }{\eta_s} 
        &\le \frac{\ep_0(1-\beta_1^s)}{C_0} + \frac{\sqrt{2}(1-\beta_1^s)}{C_0\sqrt{1-\beta_2}}\sqrt{\|\boldsymbol{\sigma}_0\|^2 + (1+\|\boldsymbol{\sigma}_1\|_{\infty}^2)G } \\
        &\le \frac{\ep_0}{C_0} + \frac{\sqrt{2} }{C_0\sqrt{1-\beta_2}}\sqrt{\|\boldsymbol{\sigma}_0\|^2 + (1+\|\boldsymbol{\sigma}_1\|_{\infty}^2)G }.
    \end{align*}
    Thus, combining with \eqref{eq:final_middle_5},
    \begin{align*}
      G &\ge L\sum_{s=1}^T \|\bar{g}_s\|^2 / \left[ \frac{\ep_0}{C_0} + \frac{\sqrt{2}}{C_0\sqrt{1-\beta_2}}\sqrt{\|\boldsymbol{\sigma}_0\|^2 + (1+\|\boldsymbol{\sigma}_1\|_{\infty}^2)G} \right].
    \end{align*}
    Hence, dividing $T$ on both sides, we have with probability at least $1-\delta$,
    \begin{align*}
        \frac{1}{T}\sum_{s=1}^T \|\bar{g}_s\|^2 \le  \frac{G}{T} \cdot \left[ \frac{\ep_0}{LC_0} + \frac{\sqrt{2}}{LC_0\sqrt{1-\beta_2}}\sqrt{\|\boldsymbol{\sigma}_0\|^2 + (1+\|\boldsymbol{\sigma}_1\|_{\infty}^2)G } \right]. 
    \end{align*}
    \hfill \BlackBox
\subsection{Proof of Corollary \ref{coro:Adam_scalar_1}}
We first provide a lemma to handle the term $\log\left(1/\beta_2^T\right)$ in $G$.
\begin{lemma}\label{lem:logbeta_t}
    Given $T \ge 2$. Suppose that $\beta_2 = 1-1/T$, then it holds $\log\left(1/\beta_2^T\right) \le 2$.
\end{lemma}
\paragraph{Proof of Lemma \ref{lem:logbeta_t}}Since $\beta_2 \in (0,1)$, we have
    \begin{align*}
        -\log \beta_2 = \log\left( \frac{1}{\beta_2} \right) \le \frac{1-\beta_2}{\beta_2} = \frac{1/T}{1-1/T} \le \frac{2}{T}.
    \end{align*}
    where we apply $\log(1/a)\le (1-a)/a,\forall a \in (0,1)$.
    With both sides multiplying $T$ we obtain the desired result. \hfill \BlackBox

\paragraph{Proof of Corollary \ref{coro:Adam_scalar_1}}
    First applying Lemma \ref{lem:logbeta_t}, we have when $\beta_2 = 1-1/T$, 
    \begin{align}
        &\log \left( \frac{\mathcal{F}(T)}{\beta_2^T} \right) \le \log \left(\mathcal{F}(T)\right) + 2 \le  \log\left(\mathrm{e}^2\mathcal{F}(T)\right). \label{eq:logbeta2_T}
    \end{align}
    Moreover, we have $1/\beta_2 \le 2$ and $1+\beta_2 \le 2$. Hence,
    \begin{align*}
        \Sigma_{\max}:= \max\left\{1, \sqrt{\frac{1+\beta_2}{\beta_2}} - 1\right\} \le \max\left\{1, \sqrt{4} - 1\right\}  \le 1.
    \end{align*}
    Hence, let us denote $G_1$ the value of $G$ in \eqref{eq:G_define} with $\beta_2 = 1-1/T$. Then combining with \eqref{eq:logbeta2_T} and $1/\beta_2 \le 2$, 
    it's then easy to verify that $G_1 \sim \mathcal{O}\left( d^2 {\rm poly}\left( \log \frac{dT}{\delta} \right)\right)$. Then, recalling the main result in Theorem \ref{thm:Adam_scalar_1}, we obtain the desired result by replacing $G$ with $G_1$.
    \begin{align*}
         \frac{1}{T}\sum_{s=1}^T \|\bar{g}_s\|^2 \le   \frac{1}{T}\cdot \frac{G_1\ep_0}{LC_0} + \frac{1}{\sqrt{T}}\cdot \frac{\sqrt{2}G_1}{LC_0}\sqrt{\|\boldsymbol{\sigma}_0\|^2 + (1+\|\boldsymbol{\sigma}_1\|_{\infty}^2)G_1 }.
    \end{align*}
    \hfill \BlackBox
\subsection{Proof of Remark \ref{coro:beta1}}
In this section we briefly discuss Remark \ref{coro:beta1}. 
We could assume that $\beta_1/\beta_2 \approx \beta_1$ based on the the experimental setting in \citep{kingma2014adam}. If $C_0 = (1-\beta_1)^a$ for some constant $a > 0$, then from Theorem \ref{thm:Adam_scalar_1} we first have
    \begin{align*}
        G \sim \mathcal{O}\left( f(x_1)-f^* + \frac{1}{(1-\beta_1)^{6-2a}} \right).
    \end{align*}
    The convergence rate in Theorem \ref{thm:Adam_scalar_1} shows that the order is dominated by $\left(G\sqrt{G}\right)/C_0$. Then we will have two cases. When $ a > 3$, $G\sqrt{G} \sim \mathcal{O}(1)$. Hence,
    \begin{align*}
        \frac{1}{T}\sum_{s=1}^T \|\nabla f(x_s)\|^2 \leq \mathcal{O}\left( \frac{1}{C_0}\right) \sim \mathcal{O}\left( \frac{1}{(1-\beta_1)^a}\right).
    \end{align*}
    The order is thus no less than $1/(1-\beta_1)^{-3}$. When $ a \le 3$, we therefore obtain that
    \begin{align*}
        \frac{1}{T}\sum_{s=1}^T \|\nabla f(x_s)\|^2 \leq \mathcal{O}\left( \frac{1}{(1-\beta_1)^3} + \frac{1}{(1-\beta_1)^a}+ \frac{1}{(1-\beta_1)^{6-a}} + \frac{1}{(1-\beta_1)^{9-2a}} \right).
    \end{align*}
    It thus leads to the minimum order of $1/(1-\beta_1)^{3}$ when $a = 3$. \hfill \BlackBox

\section{Adaptivity to the noise level of a variant of Adam
}\label{sec:improve}
Compared to the rate of the form $\mathcal{O}\left(1/T + \sigma_0/\sqrt{T}\right)$ for vanilla SGD \citep{ghadimi2013stochastic} and AdaGrad \citep{kavis2022high,attia2023sgd,alina2023high} in non-convex smooth case, although in Theorem \ref{thm:Adam_scalar_1} we establish the high probability convergence rate with the same order, the rate is not adaptive to the noise level meaning that it could not accelerate to the fast rate of $\tilde{\mathcal{O}}(1/T)$ when the noise is sufficiently low. To obtain a noise adaptation rate, we simplify Adam by dropping its corrective term for gradient square ($v_{s,i}$) in Algorithm \ref{alg:Adam_scalar}.

\begin{algorithm}\label{alg:Adam_simple}
	\caption{A variant of Adam}
	\KwIn{Horizon $T$, $x_1 \in \mathbb{R}^d$, $\beta_1, \beta_2 \in (0,1)$, $m_0 = v_0 = 0$, $\epsilon,\eta > 0$,}
	\For{$s=1,\cdots,T$}{
	Generate $g_s = (g_{s,i})_i = g(x_s) $\;
        \For{$i=1,\cdots, d$}{
        $m_{s,i}=\beta_1 m_{s-1,i} + (1-\beta_1) g_{s,i}$\;
        $v_{s,i} = \beta_2 v_{s-1,i} + (1-\beta_2)g_{s,i}^2 $\;
        $x_{s+1,i} = x_{s,i} - \frac{\eta_s}{\sqrt{v_{s,i}}+\epsilon} \cdot m_{s,i}, \tilde{\eta}_s = \eta/(1-\beta_1^s)$ \;
        }
        }
\end{algorithm}

We also include the corrective term for $\msi$ into $\tilde{\eta}_s$ and obtain that $\tilde{\eta}_s$ also satisfies $\tilde{\eta}_s \le \frac{\eta}{1-\beta_1}$.
Then we have the following convergence result.
\begin{theorem}\label{thm:Adam_simple}
    Given $T \ge 1$. Suppose that $\{x_s\}_{s \in [T]}$ is a sequence generated by Algorithm \ref{alg:Adam_simple}. Under the same conditions and parameter settings in Theorem \ref{thm:Adam_scalar_1}, then for any given $\delta \in (0,1)$, with probability at least $1-\delta$,
    \begin{align*}
        \frac{1}{T}\sum_{s=1}^T \|\nabla f(x_s)\|^2  \leq \mathcal{O}\left\{\frac{G}{TC_0}\left(\ep_0 + \frac{G(1+\|\boldsymbol{\sigma}_1\|_{\infty}^2)}{C_0} \right) + \frac{1}{\sqrt{T}}\cdot \frac{G\|\boldsymbol{\sigma}_0\|}{C_0}\right\}, 
    \end{align*}
    where $G$ follows the definition in \eqref{eq:G_define}.
\end{theorem}
\begin{remark}
    Similarly, when setting $\beta_2 = 1-1/T$, $G$ becomes $G_1$ defined in Corollary \ref{coro:Adam_scalar_1} and the convergence rate satisfies $\tilde{\mathcal{O}}\left( 1/T + \|\boldsymbol{\sigma}_0\| /\sqrt{T} \right)$ that is adaptive to the noise level. 
\end{remark}
\begin{remark}
    Our analysis could be easily extended to some other simplified version of Adam, including Adam with only corrective term for $v_{s,i}$ \citep{defossez2020simple} and Adam without corrective terms \citep{zhang2022adam}. Due to the limited spaces, we omit the details.
\end{remark}

\paragraph{Proof of Theorem \ref{thm:Adam_simple}}
    Since most of the analysis is similar to Section \ref{sec:proof}, we will briefly state the common parts and present the difference in detail. First, it's easy to verify that lemmas in Section \ref{sec:alg_lemma} still hold for Algorithm \ref{alg:Adam_simple}. We refer readers to Remark \ref{rem:lem_order_T_bound_gradient}, Remark \ref{rem:four_summation} and Remark \ref{rem:gap} (from the coming appendix) for more discussions.
    
    We also start by establishing a similar result to Proposition \ref{pro:delta_s_1}
    . First we follow the decomposition of \eqref{eq:A+B+C} with $\eta_s$ replaced by $\tilde{\eta}_s$, 
    \begin{align}
    f(y_{s+1})
   &\le f(y_s) + \sum_{i=1}^d -\tilde{\eta}_s \cdot \nabla f(y_s)_i \cdot \frac{\gsi}{\bsi} + \sum_{i=1}^d \frac{\beta_1}{1-\beta_1}\cdot \nabla f(y_s)_i \cdot \left(\frac{\tilde{\eta}_s b_{s-1,i}}{\tilde{\eta}_{s-1} b_{s,i}} -1 \right)\cdot(x_{s,i} - x_{s-1,i}) \nonumber \\
   &\quad+ \frac{L}{2}\sum_{i=1}^d \left(\tilde{\eta}_s \cdot \frac{\gsi}{\bsi}  -  \frac{\beta_1}{1-\beta_1}\cdot \left(  \frac{\tilde{\eta}_{s} b_{s-1,i}}{\tilde{\eta}_{s-1} \bsi} -1 \right)\cdot(\xsi - x_{s-1,i}) \right)^2. 
\end{align}
    Note that since $\tilde{\eta}_s$ also satisfies \eqref{eq:eta_s}, then estimations are unchanged for \textbf{A.1} in \eqref{eq:A.1}, \textbf{A.2} in \eqref{eq:A.2}, \textbf{B.1.1} in \eqref{eq:B.1.1}, \textbf{B.1.2} in \eqref{eq:B.1.2} and \textbf{C} in \eqref{eq:C}. Note that \textbf{B.1.3.2} in \eqref{eq:B.1.3_decomp} vanishes since there is no corrective term for $\vsi$. \textbf{B.1.3} becomes 
    \begin{align*}
        \textbf{B.1.3} = \sum_{i=1}^d \frac{\eta \beta_1}{1-\beta_1} \cdot \left|\left(\frac{1}{1-\beta_1^{s-1}}-\frac{1}{1-\beta_1^{s}}\right) \cdot \bgsi \cdot \frac{m_{s-1,i}}{b_{s-1,i}} \right|.
    \end{align*}
    Then we also obtain a similar result to \eqref{eq:B.1.3.1} and \eqref{eq:sum_B.1.3.1} as
    \begin{align*}
    \sum_{s=1}^t \textbf{B.1.3.1}  
    &\le \frac{d\eta \beta_1}{1-\beta_1}\sqrt{G_t}\cdot \sqrt{\frac{(1-\beta_1)(1-\beta_1^{s-1})}{(1-\beta_2)(1-\beta_1/\beta_2)}}\cdot \sum_{s=1}^t\left(\frac{1}{1-\beta_1^{s-1}}-\frac{1}{1- \beta_1^{s}}\right) \le  e\sqrt{G_t}.
    \end{align*}
    In terms of \textbf{B.2}, we see that when $\frac{\tilde{\eta}_s b_{s-1,i}}{\tilde{\eta}_{s-1} b_{s,i}} \ge 1$, it's easy to verify that $\tilde{\eta}_s/\tilde{\eta}_{s-1} = (1-\beta_1^{s-1})/(1-\beta_1^s) \le 1$.
    Combining with \eqref{eq:Sigma_max_2} and \eqref{eq:Sigma_max}, we then have
    \begin{align*}
        \left|\frac{\tilde{\eta}_s b_{s-1,i}}{\tilde{\eta}_{s-1} b_{s,i}} \right| \le \max\left\{1,\frac{1}{\sqrt{\beta_2}}-1 \right\} \le \Sigma_{\max}.
    \end{align*}
    Above all, the only difference is that \textbf{B.1.3.2} vanishes. Hence, we could still apply the induction argument and obtain the same result to Proposition \ref{pro:delta_s_1} that with probability at least $1-\delta$, for all $t\in [T],$
    \begin{align*}
        \|\bar{g}_{T+1}\|^2 \le G- L\sum_{s=1}^T \sum_{i=1}^d \tilde{\eta}_s \cdot \frac{\bgsi^2}{\asi} .
    \end{align*}
    Moreover, we have the same analysis to \eqref{eq:final_middle_5},
    \begin{align}
        G  \ge  L \sum_{s=1}^T \sum_{i=1}^d \tilde{\eta}_s \cdot \frac{\bgsi^2}{\ep_s + \sqrt{\|\tilde{\mathbf{v}}_{s}\|_1} } =L \sum_{s=1}^T \tilde{\eta}_s \cdot \frac{\|\bar{g}_s\|^2}{\ep_s + \sqrt{\|\tilde{\mathbf{v}}_{s}\|_1}} . \label{eq:final_middle_6}
    \end{align}
    We also establish an upper bound for $\|\tilde{\mathbf{v}}_{s}\|_1$ that is a little different to \eqref{eq:bound_vs}. Using basic inequality and Assumption (A3),
    \begin{align*}
        \sqrt{\|\tilde{\mathbf{v}}_{s}\|_1 } 
        &= \sqrt{\sum_{i=1}^d |\tilde{v}_{s,i}| } =  \sqrt{\sum_{i=1}^d  (1-\beta_2) \left[\sum_{j=1}^{s-1} \beta_2^{s - j} g_{j,i}^2+ (\sigma_{0,i}^2 + (1+\sigma_{1,i}^2)\bar{g}_{s,i}^2)\right]  } \nonumber \\
        &\le \sqrt{2(1-\beta_2)\sum_{i=1}^d  \sum_{j=1}^{s}\beta_2^{s - j} (\sigma_{0,i}^2 + (1+\sigma_{1,i}^2) \bar{g}_{j,i}^2) } \nonumber\\
        &= \sqrt{2(1-\beta_2)\left(\|\boldsymbol{\sigma}_0\|^2\cdot s + (1+\|\boldsymbol{\sigma}_1\|_{\infty}^2)  \sum_{j=1}^{s} \|\bar{g}_j\|^2 \right) }.
    \end{align*}
    Recalling the parameter setting in \eqref{eq:parameter_setting} where $\tilde{\eta}_s = \eta/(1-\beta_1^s) = C_0 \sqrt{1-\beta_2}/(1-\beta_1^s)$ and $\ep = \ep_0\sqrt{1-\beta_2}$, then we have for any $s \in [T]$,
    \begin{align*}
        \frac{\ep + \sqrt{\|\tilde{\mathbf{v}}_{s}\|_1 } }{\tilde{\eta}_s} 
        &\le \frac{\ep_0(1-\beta_1^s)}{C_0} +\frac{\sqrt{2}(1-\beta_1^s)}{C_0}\sqrt{\|\boldsymbol{\sigma}_0\|^2\cdot s + (1+\|\boldsymbol{\sigma}_1\|_{\infty}^2)  \sum_{j=1}^{s} \|\bar{g}_j\|^2  } \\
        &\le \frac{\ep_0}{C_0} +\frac{\sqrt{2}}{C_0}\left(\|\boldsymbol{\sigma}_0\| \cdot \sqrt{T} + \sqrt{(1+\|\boldsymbol{\sigma}_1\|_{\infty}^2)  \sum_{j=1}^{T} \|\bar{g}_j\|^2 }  \right).
    \end{align*}
    Combining with \eqref{eq:final_middle_6} and Young's inequality, we then have
    \begin{align*}
        \sum_{s=1}^T \|\bar{g}_s\|^2 
        &\le \frac{G}{L}\left[ \frac{\ep_0}{C_0} +\frac{\sqrt{2}\|\boldsymbol{\sigma}_0\|}{C_0} \cdot \sqrt{T} + \frac{\sqrt{2}}{C_0}\sqrt{(1+\|\boldsymbol{\sigma}_1\|_{\infty}^2)  \sum_{j=1}^{T} \|\bar{g}_j\|^2 }  \right] \\
        &\le \frac{G}{L}\left[\frac{\ep_0}{C_0} +\frac{\sqrt{2}\|\boldsymbol{\sigma}_0\|}{C_0}\cdot \sqrt{T} \right] + \frac{G^2(1+\|\boldsymbol{\sigma}_1\|_{\infty}^2)}{L^2C_0^2} +  \frac{1}{2} \sum_{j=1}^{T} \|\bar{g}_j\|^2 .
    \end{align*}
    Hence, re-arranging the order and dividing $T$ on both sides, it holds that with probability at least $1-\delta$,
    \begin{align*}
        \frac{1}{T}\sum_{s=1}^T \|\bar{g}_s\|^2  \le \frac{2}{T}\left(\frac{G\ep_0}{LC_0} + \frac{G^2(1+\|\boldsymbol{\sigma}_1\|_{\infty}^2)}{L^2C_0^2} \right) + \frac{1}{\sqrt{T}}\cdot \frac{2\sqrt{2}G\|\boldsymbol{\sigma}_0\|}{LC_0}.
    \end{align*}
    \hfill \BlackBox
\section*{Acknowledgments}
{This work was supported in part
	by the National Key Research and Development Program of China under grant number 2021YFA1003500, and
	NSFC under grant numbers 11971427. The corresponding author is Junhong Lin.} 

\bibliography{ref}
\appendix 
\section{Omitted proof of Section \ref{sec:tech_lemma}}
We provide detailed proofs for some of the lemmas in Section \ref{sec:tech_lemma}.
\subsection{Proof of Lemma \ref{lem:logT_2}}
Before proving the result, we introduce a lemma from \citep{defossez2020simple}.
\begin{lemma}\label{lem:logT_1}
    \citep{defossez2020simple} Given $T \ge 1$. Suppose $\{\alpha_s \}_{s \in [T]}$ is a non-negative sequence. Given $\beta_2 \in (0,1]$ and $\varepsilon > 0$, we define $\theta_s = \sum_{j=1}^s \beta_2^{s-j} \alpha_j$, then for any $t \in [T]$,
    \begin{align*}
        \sum_{s=1}^t \frac{\alpha_j}{\varepsilon + \theta_j} \le \log \left(1 + \frac{\theta_t}{\varepsilon} \right) - t \log \beta_2.
    \end{align*}
\end{lemma}
	\paragraph{Proof of Lemma \ref{lem:logT_2}}
    The first result is from \citep[Lemma A.2]{defossez2020simple}. We now prove the second result. Let $\hat{M} = \sum_{j=1}^{s} \beta_1^{s-j}$. Applying Jensen's inequality,
    \begin{align}
        \left(\sum_{j=1}^{s} \beta_1^{s-j} \alpha_j \right)^2 = \left(\hat{M} \cdot \sum_{j=1}^{s} \frac{\beta_1^{s-j}}{\hat{M}} \alpha_j \right)^2 \le \hat{M}^2 \cdot \sum_{j=1}^{s-1}\frac{\beta_1^{s-j}}{\hat{M}}  \alpha_j^2 = \hat{M} \cdot \sum_{j=1}^{s} \beta_1^{s-j} \alpha_j^2. \label{eq:Jensen}
    \end{align}
    Hence, combining with $\hat{M}=(1-\beta_1^{s})/(1-\beta_1)$ we have
    \begin{align*}
        \frac{\gamma_s^2}{\varepsilon + \theta_s} \le \frac{\hat{M}}{(1-\beta_1^s)^2} \sum_{j=1}^{s} \beta_1^{s-j} \frac{\alpha_j^2}{\varepsilon + \theta_s} 
        = \frac{1}{(1-\beta_1)(1-\beta_1^s)} \sum_{j=1}^{s} \beta_1^{s-j} \frac{\alpha_j^2}{\varepsilon + \theta_s}.
    \end{align*}
    Recalling the definition of $\theta_s$, we have $\varepsilon + \theta_s \ge \varepsilon + \beta_2^{s-j} \theta_j \ge \beta_2^{s-j}(\varepsilon+\theta_j) $. Hence, combining with $1-\beta_1 < 1-\beta_1^s$,
    \begin{align*}
        \frac{\gamma_s^2}{\varepsilon + \theta_s} \le \frac{1}{(1-\beta_1)(1-\beta_1^s)} \sum_{j=1}^{s} \left(\frac{\beta_1}{\beta_2}\right)^{s-j} \frac{\alpha_j^2}{\varepsilon + \theta_j} \le \frac{1}{(1-\beta_1)^2} \sum_{j=1}^{s} \left(\frac{\beta_1}{\beta_2}\right)^{s-j} \frac{\alpha_j^2}{\varepsilon + \theta_j}.
    \end{align*}
    Summing over $s \in [t]$, and noting that $\beta_1 < \beta_2$,
    \begin{align*}
        \sum_{s=1}^t \frac{\gamma_s^2}{\varepsilon + \theta_s} 
        &\le \frac{1}{(1-\beta_1)^2} \sum_{s=1}^t \sum_{j=1}^{s} \left(\frac{\beta_1}{\beta_2}\right)^{s-j} \frac{\alpha_j^2}{\varepsilon + \theta_j} = \frac{1}{(1-\beta_1)^2} \sum_{j=1}^t \frac{\alpha_j^2}{\varepsilon + \theta_j} \sum_{s=j}^{t} \left(\frac{\beta_1}{\beta_2}\right)^{s-j}   \\
        &\le \frac{1}{(1-\beta_1)^2(1-\beta_1/\beta_2)} \sum_{j=1}^t \frac{\alpha_j^2}{\varepsilon + \theta_j}.
    \end{align*}
    Finally applying Lemma \ref{lem:logT_1}, we obtain the desired result. \hfill \BlackBox

\section{Omitted proof of Section \ref{sec:alg_lemma}}
We provide detailed proofs for some of the lemmas in Section \ref{sec:alg_lemma}.
\subsection{Proof of Lemma \ref{lem:OrderT_bound_gradient}}
    Noting that from the smoothness of $f$,
    \begin{align}
        \|\nabla f(x_s)\| 
        &\le \|\nabla f(x_{s-1}) \| + \|\nabla f(x_s) - \nabla f(x_{s-1})\| \le  \|\nabla f(x_{s-1}) \| + L \|x_s - x_{s-1}\|. \label{eq:update_1}
    \end{align}
    Recalling the updated rule, we have for any $i \in [d]$,
    \begin{align*}
        |x_{s,i} - x_{s-1,i}| = \left|\frac{\eta_{s-1}}{\sqrt{v_{s-1,i}}+\ep_{s-1}} m_{s-1,i}  \right| \le \eta_{s-1} \left| \frac{m_{s-1,i}}{\sqrt{v_{s-1,i}}} \right|.
    \end{align*}
    Recalling the definitions of $m_{s-1,i}$ and $v_{s-1,i}$ in Algorithm \ref{alg:Adam_scalar}. Then applying \eqref{eq:Jensen} with $\hat{M}$ replaced by $\tilde{M}=\sum_{j=1}^{s-1} \beta_1^{s-1-j}$,
    \begin{align*}
         \left| \frac{m_{s-1,i}}{\sqrt{v_{s-1,i}}} \right|  
         &= \sqrt{\frac{(1-\beta_1)^2\left(\sum_{j=1}^{s-1} \beta_1^{s-1-j} g_{j,i} \right)^2 }{(1-\beta_2)\sum_{j=1}^{s-1} \beta_2^{s-1-j} g^2_{j,i}} }  \le \frac{1-\beta_1}{\sqrt{1-\beta_2}}\sqrt{\tilde{M}\cdot \frac{ \sum_{j=1}^{s-1} \beta_1^{s-1-j} g_{j,i}^2 }{\sum_{j=1}^{s-1} \beta_2^{s-1-j} g^2_{j,i}} } \nonumber \\
        &\le  \frac{1-\beta_1}{\sqrt{1-\beta_2}}\sqrt{\tilde{M}\cdot \sum_{j=1}^{s-1} \left(\frac{\beta_1}{\beta_2}\right)^{s-1-j} } =  \frac{1-\beta_1}{\sqrt{1-\beta_2}} \sqrt{\frac{1-\beta_1^{s-1}}{1-\beta_1} \cdot \frac{1-(\beta_1/\beta_2)^{s-1}}{1-\beta_1/\beta_2}} \nonumber \\
        &\le  \sqrt{\frac{(1-\beta_1)(1-\beta_1^{s-1})}{(1-\beta_2)(1-\beta_1/\beta_2)} },
    \end{align*}
    where the last inequality applies $\beta_1 < \beta_2$. Then combining with the parameter setting where $\eta_{s-1} = \eta\sqrt{1-\beta_2^{s-1}}/(1-\beta_1^{s-1})$ and $\eta=C_0\sqrt{1-\beta_2}$, we further have
    \begin{align}
        |x_{s,i} - x_{s-1,i}| 
        &\le \eta \sqrt{\frac{(1-\beta_1)(1-\beta_2^{s-1})}{(1-\beta_2)(1-\beta_1/\beta_2)(1-\beta_1^{s-1})} } \le C_0 \sqrt{\frac{1}{1-\beta_1/\beta_2}}. \label{eq:bound_ms-1/bs-1_middle}
    \end{align}
    Hence, we have $\|x_{s,i} - x_{s-1,i}\|_{\infty} \le C_{0} \sqrt{\frac{1}{1-\beta_1/\beta_2} }$. Combining with $\|x_{s}-x_{s-1}\| \le \sqrt{d}\|x_{s,i} - x_{s-1,i}\|_{\infty}$ and \eqref{eq:update_1}, 
    \begin{align*}
        \|\nabla f(x_s)\|_{\infty} &\le \|\nabla f(x_s)\| \le \|\nabla f(x_{s-1}) \| + LC_{0}\sqrt{\frac{d}{1-\beta_1/\beta_2} } \\
        &\le \|\nabla f(x_1)\| + LC_{0}\sqrt{\frac{d}{1-\beta_1/\beta_2} }  \cdot s. 
    \end{align*}
    \hfill \BlackBox
\begin{remark}\label{rem:lem_order_T_bound_gradient}
    Note that if we take the step-size $\tilde{\eta}_s$ in Algorithm \ref{alg:Adam_simple}, the estimation of \eqref{eq:bound_ms-1/bs-1_middle} remains unchanged. Hence, Lemma \ref{lem:OrderT_bound_gradient} also holds for Algorithm \ref{alg:Adam_simple}.
\end{remark}

\subsection{Proof of Lemma \ref{lem:gradient_xs_ys}}
    Applying the norm inequality and the smoothness of $f$, 
    \begin{align*}
        \|\nabla f(x_s) \| &\le \|\nabla f(y_s) \| + \|\nabla f(x_s) - \nabla f(y_s) \|  \le \|\nabla f(y_s) \| + L\|y_s-x_s\|.
    \end{align*}
    Using \eqref{eq:define_y_s} and \eqref{eq:bound_ms-1/bs-1_middle},
    \begin{align*}
        \|\nabla f(x_s) \| \le \|\nabla f(y_s) \| + \frac{L\beta_1}{1-\beta_1}\|x_s - x_{s-1}\| \le \|\nabla f(y_s) \| + \frac{L\beta_1}{1-\beta_1} \cdot C_{0}\sqrt{\frac{d}{1-\beta_1/\beta_2} }.
    \end{align*}
    \hfill \BlackBox

\subsection{Proof of Lemma \ref{lem:sum_1}}
     Recalling the updated rule and the definition of $\bsi$ in \eqref{eq:bsi}, using $\ep_s^2 = \ep^2(1-\beta_2^s) \ge \ep^2(1-\beta_2)$, 
    \begin{align}
        b_{s,i}^2 \geq  v_{s,i}^2 + \ep_s^2 \geq (1-\beta_2) \left( \sum_{j=1}^s \beta_2^{s-j}g_{j,i}^2 + \ep^2\right), \quad \text{and} \quad m_{s,i} = \left(1-\beta_1\right)\sum_{j=1}^s \beta_1^{s-j} g_{j,i}. \label{eq:express_bsi}
    \end{align}
    Also, applying the basic inequality, $\beta_2\leq 1$ and Assumption (A3),
    \begin{align}
   \sum_{s=1}^t \beta_2^{t-s}\gsi^2 
    \le & 2\sum_{s=1}^t \left(\bgsi^2 + \xisi^2 \right) \le  2 \sum_{s=1}^t \left(\sigma_{0,i}^2 + (1+\sigma_{1,i}^2)\bgsi^2 \right) \nonumber \\
    \le & 2 \left(  \|\boldsymbol{\sigma}_0\|_{\infty}^2 t
    +
   (1+\|\boldsymbol{\sigma}_{1}\|_{\infty}^2)  \sum_{s=1}^t \|\bar{g}_s\|_\infty^2  \right). \nonumber
    \end{align}
    Combining with Lemma   \ref{lem:OrderT_bound_gradient},
    \begin{equation}\label{eq:log_F(t)_1}
    1+ \frac{1}{\ep^2} \sum_{s=1}^t \beta_2^{t-s}\gsi^2   \leq   \underbrace{ 1+ \frac{2}{\ep^2}\left[(\|\boldsymbol{\sigma}_0 \|_{\infty}^2 + 2(1+\|\boldsymbol{\sigma}_1\|_{\infty}^2) \|\bar{g}_1\|^2) \cdot t + \frac{2(1+\|\boldsymbol{\sigma}_1 \|_{\infty}^2)dL^2C_0^2 }{1-\beta_1/\beta_2} \cdot t^3 \right]}_{\mathcal{F}(t)}.
    \end{equation}
    \paragraph{Proof of the first summation} Using \eqref{eq:express_bsi}, for any $i \in [d]$,
    \begin{align*}
        \sum_{s=1}^t\frac{\gsi^2}{\bsi^2} 
         \le \frac{1}{1-\beta_2}\sum_{s=1}^t\frac{\gsi^2}{\ep^2 + \sum_{j=1}^s \beta_2^{s-j}g_{j,i}^2}.
    \end{align*}
Applying Lemma \ref{lem:logT_1},
    \begin{align}
        \sum_{s=1}^t\frac{\gsi^2}{\bsi^2} 
        &\le \frac{1}{1-\beta_2} \left[  \log \left( 1+ \frac{1}{\ep^2} \sum_{s=1}^t \beta_2^{t-s}\gsi^2 \right) -t \log\beta_2 \right]  \label{eq:gsi/bsi}.
            \end{align}
   Introducing with \eqref{eq:log_F(t)_1},
    combining with $\beta_2 \le 1$ and noting that $\mathcal{F}(t)$ is increasing and $\beta_2^t$ is decreasing with respect to $t$, we obtain the first desired result
     by 
     summing over $i \in [d]$.
    \paragraph{Proof of the second summation } 
    Following from \eqref{eq:express_bsi},
    \begin{align*}
        \sum_{s=1}^{t-1}\frac{\msi^2}{\bsi^2} 
        \le \frac{(1-\beta_1)^2}{1-\beta_2}\cdot \sum_{s=1}^{t-1}\frac{\left(\sum_{j=1}^s \beta_1^{s-j} g_{j,i}\right)^2}{\ep^2 +  \sum_{j=1}^s \beta_2^{s-j}g_{j,i}^2}.
    \end{align*}
    Applying Lemma \ref{lem:logT_2}, and introducing with \eqref{eq:log_F(t)_1} with $t$ replaced by $t-1$, 
    \begin{align}
        \sum_{s=1}^{t-1}\frac{\msi^2}{\bsi^2}
        &\le \frac{(1-\beta_1)^2}{1-\beta_2}\cdot \frac{1}{(1-\beta_1)(1-\beta_1/\beta_2)}\left[ \log\left( 1+ \frac{1}{\ep^2} \sum_{s=1}^{t-1} \beta_2^{t-1-s}g_{s,i}^2 \right) -(t-1) \log\beta_2 \right] \nonumber
        \\
        &\le \frac{(1-\beta_1)}{(1-\beta_2)(1-\beta_1/\beta_2)} \left[ \log(\mathcal{F}(t-1)) - (t-1)\log \beta_2 \right] \nonumber. 
    \end{align}
     Noting that $\mathcal{F}(t)$ is increasing and $\beta_2^t$ is decreasing with respect to $t$, we obtain the second desired result
    by 
    summing over $i \in [d]$.
    \paragraph{Proof of the third summation} Following from \eqref{eq:express_bsi},
    \begin{align*}
        \sum_{s=1}^{t-1}\frac{\msi^2}{b^2_{s+1,i}} 
        &\le \sum_{s=1}^{t-1}\frac{\left[(1-\beta_1)\sum_{j=1}^{s} \beta_1^{s-j} g_{j,i}\right]^2}{\ep^2(1-\beta_2) + (1-\beta_2) \sum_{j=1}^{s+1} \beta_2^{s+1-j}g_{j,i}^2} \\
        &\le \sum_{s=1}^{t-1}\frac{(1-\beta_1)^2 \left(\sum_{j=1}^{s} \beta_1^{s-j} g_{j,i}\right)^2}{\ep^2(1-\beta_2) + (1-\beta_2)\beta_2 \sum_{j=1}^{s} \beta_2^{s-j}g_{j,i}^2}.
    \end{align*}
    Applying Lemma \ref{lem:logT_2},
    \begin{align}
        \sum_{s=1}^{t-1}\frac{\msi^2}{b^2_{s+1,i}} 
        &\le \frac{(1-\beta_1)^2}{(1-\beta_2)\beta_2}\cdot \sum_{s=1}^{t-1}\frac{\left(\sum_{j=1}^{s} \beta_1^{s-j} g_{j,i} \right)^2}{\frac{\ep^2}{\beta_2} + \sum_{j=1}^{s} \beta_2^{s-j}g_{j,i}^2} \nonumber \\
        &\le \frac{(1-\beta_1)^2}{(1-\beta_2)\beta_2}\cdot \frac{1}{(1-\beta_1)(1-\beta_1/\beta_2)}\left[ \log\left( 1+ \frac{\beta_2}{\ep^2} \sum_{s=1}^{t-1} \beta_2^{t-1-s}g_{s,i}^2 \right) -(t-1) \log\beta_2 \right] \label{eq:log_F(t)_4}.
    \end{align}
    Using $\beta_2 \le 1$ and then introducing with \eqref{eq:log_F(t)_1},
    \begin{align*}
        \sum_{s=1}^{t-1}\frac{\msi^2}{b^2_{s+1,i}} 
        &\le \frac{(1-\beta_1)}{\beta_2(1-\beta_2)(1-\beta_1/\beta_2)}\left[ \log(\mathcal{F}(t-1)) - (t-1)\log \beta_2 \right] \nonumber \\
        &\le \frac{(1-\beta_1)}{\beta_2(1-\beta_2)(1-\beta_1/\beta_2)}\log\left( \frac{\mathcal{F}(T)}{\beta_2^T} \right). 
    \end{align*}
    Summing over $i\in [d]$, we obtain the third desired result.
    \paragraph{Proof of the fourth summation}
     Following the definition of $\hmsi$ from \eqref{eq:hat_msi}, and combining with \eqref{eq:express_bsi},
    \begin{align*}
        \sum_{s=1}^{t-1}\frac{\hmsi^2}{\bsi^2} 
        &\le \frac{(1-\beta_1)^2}{1-\beta_2}\cdot  \sum_{s=1}^{t-1}\frac{\left( \frac{1}{1-\beta_1^s} \sum_{j=1}^s \beta_1^{s-j} g_{j,i}\right)^2}{\ep^2 + \sum_{j=1}^s \beta_2^{s-j}g_{j,i}^2}.
    \end{align*}
    Applying Lemma \ref{lem:logT_2}, and using \eqref{eq:log_F(t)_1}, 
    \begin{align}
        \sum_{s=1}^{t-1}\frac{\hmsi^2}{\bsi^2}
        &\le \frac{(1-\beta_1)^2}{1-\beta_2}\cdot \frac{1}{(1-\beta_1)^2(1-\beta_1/\beta_2)}\left[ \log\left( 1+ \frac{1}{\ep^2} \sum_{s=1}^{t-1} \beta_2^{t-1-s}g_{s,i}^2 \right) -(t-1) \log\beta_2 \right] \label{eq:log_F(t)_5}\\
        &\le \frac{1}{(1-\beta_2)(1-\beta_1/\beta_2)} \left[ \log(\mathcal{F}(t-1)) - (t-1)\log \beta_2 \right] \nonumber\\
        &\le \frac{1}{(1-\beta_2)(1-\beta_1/\beta_2)}\log\left( \frac{\mathcal{F}(T)}{\beta_2^T} \right). \nonumber
    \end{align}
     Summing over $i\in [d]$, we obtain the fourth desired result.
    \hfill \BlackBox
\begin{remark}\label{rem:four_summation}
    If we consider Algorithm \ref{alg:Adam_simple}, we should replace $\ep_s$ with the constant $\ep$, and the estimation 
    \eqref{eq:express_bsi} still holds. 
Thus, using the same arguments, one can verify that Lemma \ref{lem:sum_1} still holds for Algorithm \ref{alg:Adam_simple}. 
\end{remark}
\subsection{Proof of Lemma \ref{lem:gap_as_bs}}
    Recalling the definition in \eqref{eq:prox_as} and Assumption (A3), we have
    \begin{align*}
        \left| \frac{1}{\asi} - \frac{1}{\bsi} \right|  
        &= \frac{\left|\sqrt{\vsi}-  \sqrt{\tilde{v}_{s,i}}\right|}{\asi \bsi} \le \frac{(1-\beta_2)}{\asi \bsi}\frac{\left|\gsi^2 - \bgsi^2  - \sigma_{0,i}^2 - \sigma_{1,i}^2 \bgsi^2 \right|}{\sqrt{v_{s,i}}+ \sqrt{\tilde{v}_{s,i}}} \nonumber \\
        &\le \frac{(1-\beta_2)}{\asi \bsi}\frac{\left|\gsi - \bgsi \right|\left|\gsi + \bgsi  \right| + \sigma_{0,i}^2 + \sigma_{1,i}^2 \bgsi^2  }{\sqrt{v_{s,i}}+ \sqrt{\tilde{v}_{s,i}}} \nonumber \\
        &\le \frac{(1-\beta_2)}{\asi \bsi}\frac{\left|\xisi\right|(\left|\gsi\right| + \left|\bgsi  \right|) + \sigma_{0,i}^2 + \sigma_{1,i}^2 \bgsi^2 }{\sqrt{\beta_2 v_{s-1,i} + (1-\beta_2)g_{s,i}^2} + \sqrt{\beta_2 v_{s-1,i} + (1-\beta_2)(\sigma_{0,i}^2 + (1+\sigma_{1,i}^2) \bgsi^2 ) }} \nonumber \\
        &\le \frac{\sqrt{1-\beta_2}\left( |\xisi| + \sqrt{\sigma_{0,i}^2 + \sigma_{1,i}^2 \bgsi^2}\right)}{\asi \bsi} \le \frac{2\sqrt{1-\beta_2}}{\asi \bsi}\sqrt{\sigma_{0,i}^2 + \sigma_{1,i}^2 \bgsi^2}.
    \end{align*}
    The second result also follows from the same analysis. We first have
    \begin{align}
        \left|\frac{1}{b_{s-1,i}} - \frac{1}{\asi} \right| &= \frac{\left|\sqrt{\tilde{v}_{s,i}} - \sqrt{v_{s-1,i}}+ (\ep_s - \ep_{s-1})\right|}{b_{s-1,i}\asi}  \nonumber\\
        &\le \frac{1}{b_{s-1,i}\asi} \frac{(1-\beta_2)\left| \sigma_{0,i}^2 + (1+\sigma_{1,i}^2)\bgsi^2 - v_{s-1,i} \right|}{\sqrt{\tilde{v}_{s,i}} + \sqrt{v_{s-1,i}}}  + \frac{\left|\ep_s - \ep_{s-1}\right|}{b_{s-1,i}\asi} \label{eq:gap_ep} \\
        &=\frac{1}{b_{s-1,i}\asi} \frac{(1-\beta_2)\left| \sigma_{0,i}^2 + (1+\sigma_{1,i}^2)\bgsi^2 - v_{s-1,i} \right|}{\sqrt{\tilde{v}_{s,i}} +\sqrt{v_{s-1,i}}} + \frac{\ep\left(\sqrt{1-\beta_2^s}-\sqrt{1-\beta_2^{s-1}}\right)}{b_{s-1,i}\asi} \nonumber \\
        &\le \frac{1}{b_{s-1,i}\asi} \frac{(1-\beta_2)\left| \sigma_{0,i}^2 + (1+\sigma_{1,i}^2)\bgsi^2 - v_{s-1,i} \right|}{\sqrt{\tilde{v}_{s,i}}+ \sqrt{v_{s-1,i}}} + \frac{\ep\sqrt{\beta_2^{s-1}(1-\beta_2)}}{b_{s-1,i}\asi}, \nonumber
    \end{align}
    where the last inequality applies that $\sqrt{a}-\sqrt{b} \le \sqrt{a-b}, \forall 0 \le b \le a$. Note that if $ \sigma_{0,i}^2 + (1+\sigma_{1,i}^2)\bgsi^2 > v_{s-1,i}$, combining with $\beta_2 \le 1$, we obtain that
    \begin{align*}
        \left|\frac{1}{b_{s-1,i}} - \frac{1}{\asi} \right|  
        &\le \frac{1}{b_{s-1,i}\asi} \cdot \frac{(1-\beta_2)  (\sigma_{0,i}^2 + (1+\sigma_{1,i}^2)\bgsi^2) }{\sqrt{\tilde{v}_{s,i}}+ \sqrt{v_{s-1,i}}} + \frac{\ep\sqrt{\beta_2^{s-1}(1-\beta_2)}}{b_{s-1,i}\asi} \\
        &\le \frac{\sqrt{1-\beta_2 } }{b_{s-1,i}\asi}\left(\sqrt{\sigma_{0,i}^2 + (1+\sigma_{1,i}^2)\bgsi^2} + \ep \right).
    \end{align*}
    If $ \sigma_{0,i}^2 + (1+\sigma_{1,i}^2)\bgsi^2 \le v_{s-1,i}$, combining with $1-\beta_2 \le \sqrt{1-\beta_2}$,
    \begin{align*}
        \left|\frac{1}{b_{s-1,i}} - \frac{1}{\asi} \right|
        &\le \frac{1}{b_{s-1,i}\asi} \cdot \frac{(1-\beta_2) v_{s-1,i} }{\sqrt{\tilde{v}_{s,i}} +  \sqrt{v_{s-1,i}}}  + \frac{\ep\sqrt{\beta_2^{s-1}(1-\beta_2)}}{b_{s-1,i}\asi}\\
        &\le  \frac{\sqrt{v_{s-1,i}}(1-\beta_2)}{b_{s-1,i}\asi}+ \frac{\ep\sqrt{1-\beta_2}}{b_{s-1,i}\asi} \le \frac{\sqrt{1-\beta_2}}{b_{s-1,i}\asi}( \sqrt{v_{s-1,i}} + \ep).
    \end{align*}
    Combining the two cases, we obtain the desired result.\hfill\BlackBox
\begin{remark}\label{rem:gap}
    Noting that when performing Algorithm \ref{alg:Adam_simple}, the estimation for the gap $\left|1/\asi - 1/\bsi \right|$ remains unchanged. Since $\ep_s$ is a constant in Algorithm \ref{alg:Adam_simple}, the second term in \eqref{eq:gap_ep} vanishes. Hence, the gap $\left|1/b_{s-1,i} - 1/\asi \right|$ becomes smaller but we could still use the estimation in Lemma \ref{lem:gap_as_bs} as an upper bound. 
\end{remark}

\end{document}